\documentclass[11pt]{amsart}
\usepackage{graphicx}

\usepackage[utf8]{inputenc}		% Codificacao do documento (conversão automática dos acentos)
\usepackage{indentfirst}		% Indenta o primeiro parágrafo de cada seção.
\usepackage{graphicx}			% Inclusão de gráficos
 \usepackage{amsbsy}			% para símbolos matemáticos em negrito
 \usepackage{amscd}			% para diagramas
 \usepackage{amsfonts}			% fontes AMS
 \usepackage{amsmath}			% facilidades matemáticas
 \usepackage{amssymb}			% para os símbolos mais antigos
 \usepackage{amstext}			% para fragmentos tipo texto em modo matemático
\usepackage{amsthm}			% para teoremas
\usepackage{hyperref}			% Amplo suporte para hipertexto em LaTeX
\usepackage{cleveref}			% Referência cruzada inteligente
\usepackage{dsfont}			% para o estilo de conjuntos de números $\mathds{R}$
 \usepackage{ifthen}			% comandos de condição em LaTeX
\usepackage{listings}           	% para inserir códigos de outras linguagens de programação
 \usepackage{lscape}             	% para imprimir alguma página no formato paisagem
\usepackage{mathabx}			% conjunto de simbolos matemáticos
 \usepackage{mathrsfs}			% suporte para fontes RSFS
 \usepackage{pdfpages}           	% para inserir páginas PDF no texto

 \usepackage{verbatim}

\usepackage{latexsym}
\usepackage{amsthm}
\usepackage{times}
\usepackage{graphicx}
\usepackage{graphics}
\usepackage{hyperref}         % para acrescentar links
\usepackage{fancybox}
\usepackage{fancyhdr}
\usepackage{ae}
\usepackage[all]{xy}

\newcommand{\R}{\mathbb{R}}

\newcommand{\lie}{\mathfrak{g}}
\newcommand{\Ad}{\text{Ad}}
\newcommand{\cqd}{\begin{flushright}$\Box$\end{flushright}}
\newcommand{\m}{\mathfrak{m}}
\newcommand{\sub}{\mathfrak{k}}
\newcommand{\h}{\mathfrak{h}}
\newcommand{\p}{\mathfrak{p}}
\newcommand{\q}{\mathfrak{q}}
\newcommand{\cvg}{\text{\bf g}_t}
\newcommand{\cvh}{\text{\bf h}_t}
\newcommand{\cvk}{\text{\bf k}_t}
\newcommand{\cvm}{\text{\bf m}_t}
\newcommand{\cvn}{\text{\bf n}_t}

\newtheorem{theorem}{Theorem}[section]
\newtheorem{proposition}[theorem]{Proposition}
\newtheorem{lemma}[theorem]{Lemma}
\newtheorem{corollary}[theorem]{Corollary}
\theoremstyle{definition}
\newtheorem{definition}[theorem]{Definition}
\newtheorem{example}[theorem]{Example}
\theoremstyle{remark}
\newtheorem{remark}[theorem]{Remark}
\numberwithin{equation}{section}

\begin{document}

\title[Yamabe Problem on Maximal Flag Manifolds]{Bifurcation and Local Rigidity of Homogeneous Solutions to the Yamabe Problem on Maximal Flag Manifolds}

 % Delete if not wanted.
\author{Lino Grama}
\address{Instituto de Matemática, Estatística e Computação Científica, Universidade Estadual de Campinas,
Campinas, SP}
\email{linograma@gmail.com}
\author{Kennerson N. S. Lima}
\address{Unidade Acadêmica de Matemática, Universidade Federal de Campina Grande, 
Paraíba, PB}
\email{kennerson@mat.ufcg.edu.br}

\begin{abstract}
In this work we determine bifurcation instants for 1-parameter families of solutions to the Yamabe problem defined on maximal flag manifolds. We also study the local rigidity points, namely, a isolated solution of the Yamabe problem. 

\end{abstract}

% keywords can be removed
\maketitle
%\tableofcontents

\section{Introduction}
Given a compact orientable Riemannian manifold $(M,g)$ with dimension $m\geq3$, the {\it Yamabe problem} concerns the existence of constant scalar curvature metrics on $M$ conformal to $g$. Solutions to this problem, called {\it Yamabe metrics}, can be characterized variationally as critical points of the Hilbert-Einstein functional 
\begin{equation*}
\mathcal{A}(g)=\dfrac{1}{\textrm{Vol}(g)}\int_M\textrm{scal}(g)\textrm{vol}_g, 
\end{equation*}
restricted to the set $[g]$ of metrics conformal to $g$. The existence of such solutions are consequence of the successive works of Yamabe \cite{yamabe}, Trudinger \cite{trundinger}, Aubin \cite{aubin} and Schoen \cite{schoen}.

The families that we study in this work are formed by homogeneous metrics which are trivial solutions to the Yamabe problem. 

A {\em bifurcation point} of a family of solutions to the Yamabe problem is an accumulation point of other solutions of the Yamabe problem conformal to the homogeneous solutions. On the other hand, a {\em local rigidity point} is an isolated solution of the Yamabe problem on its conformal class, i.e., it is not a bifurcation point. 

The bifurcation theory applied here is based on finding bifurcation instants for a given 1-parameter family of homogeneous metrics $g_t, 0<t\leq 1$. In order to do this we proceed a careful analysis of the occurrence of {\it jump} in the Morse index of the metric $g_t$. The main ingredient here is an explicity description of the Laplacian defined on $(M,g_t)$ and the scalar curvature of $g_t$.

Given a {\it Riemannian submersion with totally geodesic fibers} one can construct an 1-parameter family of  Riemannian submersions by scaling the original metric of the total space in the direction of the fibers. This family is called {\it canonical variation}. A special example of Riemannian submersion with totally geodesic fibers are the so-called {\it homogeneous fibrations}.

The study of bifurcation points and local rigidity instants for canonical variations of the round metric on spheres was provided by Bettiol and Piccione in \cite{bettiol}. These families of deformed metrics were constructed from homogeneous metrics on total spaces of Hopf fibrations (the total spaces equipped with such deformed metrics are also often referred to as {\it Berger spheres}).

In this work, we will study bifurcation points and local rigidity instants for canonical variations on maximal flag manifolds. These flag manifolds are viewed as total spaces of the following homogeneous fibrations:

\begin{itemize}
\item[(i)] $SU(n)/T^{n-1}\cdots SU(n+1)/T^{n} \rightarrow SU(n+1)/S(U(1)\times U(n))=\mathbb{CP}^{n}, n\geq 2;$
\item[(ii)] $SO(2n)/T^{n}\cdots SO(2n+1)/T^{n} \rightarrow SO(2n+1)/SO(2n)=S^{2n}, n\geq 2, n\neq 3;$
\item[(iii)] $SU(n)/T^{n-1}\cdots Sp(n)/T^{n} \rightarrow Sp(n)/U(n), n\geq 3;$
\item[(iv)]$SU(n)/T^{n-1}\cdots SO(2n)/T^{n} \rightarrow SO(2n)/U(n), n\geq 4;$
\item[(v)]$SO(4)/T\cong S^2\times S^2\cdots G_2/T \rightarrow G_2/SO(4).$
\end{itemize}  

More precisely, we equip each total space above with a {\it normal homogeneous metric} and deform these metrics by shrinking the fibers by a factor $t^2$, $0<t\leq 1,$ getting the canonical variations $$(SU(n+1)/T^{n},\cvg), (SO(2n+1)/T^{n},\cvh), (Sp(n)/T^{n},\cvk), (SO(2n)/T^{n},\cvm), (G_2/T,\cvn).$$

 Recall that a normal homogeneous metric on $G/K$ is obtained from the restriction to the tangent space (isotropy representation) of a bi-invariant inner product on the Lie algebra $\text{Lie}(G)=\lie$. The maximal flag manifolds are homogeneous spaces such that $K$ is compact, hence they admit bi-invariant metrics and hence normal homogeneous metrics.

A {\it degeneracy point} for a given canonical variation $g_t, t>0$, is a degenerate critical point $g_{t_{\ast}}$ for the Hilbert-Einstein functional, at some $t_{\ast}>0$. It is established that every bifurcation point is a degeneracy point for this canonical variation, however not all degeneracy point is a bifurcation point. Since the Morse index of each $\cvg, \cvh, \cvk, \cvm, \cvn$, changes as $t$ pass by a {\it degeneracy instants} $t_{\ast}\in ]0,1[$, we prove existence of new solutions to the Yamabe problem accumulating at $\text{\bf g}_{t_{\ast}}, \text{\bf h}_{t_{\ast}}, \text{\bf k}_{t_{\ast}}, \text{\bf m}_{t_{\ast}},\text{\bf n}_{t_{\ast}},$ respectively. Such instants $t_{\ast}$ are called {\it bifurcation instants} and $\text{\bf g}_{t_{\ast}}, \text{\bf h}_{t_{\ast}}, \text{\bf k}_{t_{\ast}}, \text{\bf m}_{t_{\ast}},\text{\bf n}_{t_{\ast}}$ are bifurcation points for the canonical variations $\cvg, \cvh, \cvk, \cvm, \cvn$, respectively. We also determine the local rigidity instants for each family $\cvg,\cvh,\cvk,\cvm,\cvn$ in the interval $]0,1]$, see Proposition \ref{criteriorig}.

Our results also can be understood from the view point of dynamical systems, where {\it bifurcation} means a topological or qualitative change in the structure of the set of fixed points of a 1-parameter family of systems when we vary this paremeter. Critical points of the Hilbert-Einstein functional in a conformal class $[g]$ are fixed points of the so-called {\it Yamabe flow}, the corresponding $L^2$-gradient flow of the Hilbert-Einstein functional, which gives a dynamical system in this conformal class. Hence, the bifurcation results above mentioned can be interpreted as a local change in the set of fixed points of the Yamabe flow near homogeneous metrics (which are always fixed points) when varying the conformal class $[g]$ with $g$ in one of the families $\cvg, \cvh, \cvk, \cvm, \cvn$.

In order to prove the bifurcation results claimed, we analyze the second variation of the Hilbert-Einstein functional at every homogeneous metrics $\cvg, \cvh, \cvk, \cvm, \cvn$, defined on the total spaces in (i)-(v). Given the second variation formula for this functional, such analysis is equivalent to compare the eigenvalues $\lambda^{kj}(t)$ (see expression \eqref{kj}, Corollary \ref{deltat}) of the Laplacian $\Delta_t$ of $\cvg, \cvh, \cvk, \cvm, \cvn$, respectively, with the scalar curvatures $\text{scal}(\cvg),\text{scal}(\cvh),\text{scal}(\cvk),\text{scal}(\cvm),\text{scal}(\cvn)$. More precisely, a critical point $\cvg, \cvh, \cvk, \cvm$ or $\cvn$ is degenerate if and only if $\text{scal}(t)\neq 0$ and $\dfrac{\text{scal}(t)}{m-1}$ is an eigenvalue of $\Delta_t.$ 

The spectrum of the Laplacian $\Delta_t$ of the canonical variations of Riemannian submersions with totally geodesic fibers is well-understood. Roughly speaking, it consists of linear combinations (that depends on $t$) of eigenvalues of the original metric on the total space with eigenvalues of the fibers. In particular, we compute a formula for the first positive eigenvalues of $\Delta_{\cvg}$ and $\Delta_{\cvh}$, and a lower and upper bounds for the first positive eigenvalues of $\Delta_{\cvk}$, $\Delta_{\cvm}$ and $\Delta_{\cvn}$, see (Proposition \ref{primeiroaut}). Combining this knowledge of the spectra of $\Delta_{\cvg}$, $\Delta_{\cvh}$, $\Delta_{\cvk}$, $\Delta_{\cvm}$ and $\Delta_{\cvn}$ with the formula for the scalar curvatures $\text{scal}(\cvg),\text{scal}(\cvh),\text{scal}(\cvk),\text{scal}(\cvm),\text{scal}(\cvn)$ we are able to identify all degeneracy instants for $\cvg, \cvh, \cvk, \cvm, \cvn$ and prove existence of bifurcation at all degeneracy instants in the interval $]0,1[$. The local rigidity instants for these canonical variations are also determined for all other $0<t<1$. We also compute the Morse index of $\cvg$ and $\cvh$, for  all $0<t\leq 1$.

With respect to the spectra of maximal flag manifolds endowed with a normal homogeneous metric we use the description of Yamaguchi (see \cite{yamaguchi}). In the case of the spectrum of the Laplacian on compact isotropy irreducible Hermitian symmetric space $G/H$ we have used in this work the description presented by Urakawa and Sugiura (see  \cite{Urakawa} and \cite{sugiura}). It is important to point out that the homogeneous metrics $\cvg, \cvh, \cvk, \cvm, \cvn$ are {\em not} normal metrics on the respectives maximal flag manifolds for $0<t<1$ and the spectra of their Laplacian $\Delta_t$ are well understood in general.  
 
Exploring the existence of infinitely many bifurcations of $\cvg, \cvh, \cvk, \cvm, \cvn$ in Section \ref{Section3.3} we obtain for each family of homogeneous spaces the existence of a subset $\mathcal{G}\subset ]0,1[$, accumulating at $0$ and such that for each $t\in\mathcal{G}$, there are at least $3$ solutions to the Yamabe problem in each conformal class $[\cvg], [\cvh], [\cvk], [\cvm], [\cvn].$

This paper is organized as follows. In Section \ref{Section1} we introduce the variational characterization of the Yamabe problem and the basic notions of bifurcation and local rigidity. In Section \ref{flag} is established the definitions of generalized flag manifold, invariant metrics and isotropy representation.  in Section \ref{cvariation} we the Laplacian of the Riemannian submersions with totally geodesic fibers, the definition of the canonical variations and the description of their spectra, as well as the construction of the canonical variations $\cvg, \cvh, \cvk, \cvm$ and $\cvn$. In Section \ref{section5} we obtain the formula for scalar curvature of $\cvg, \cvh, \cvk, \cvm$ and $\cvn$. We also prove in Section \ref{section5} the main results of this paper namely bifurcation and local rigidity of solutions of the Yamabe problem for the families $\cvg, \cvh, \cvk, \cvm$ and $\cvn$ respectively. We also determine the Morse index of $\cvg$ and $\cvh.$ The last section of this chapter contains multiplicity of solutions to the Yamabe problem for all the families above mentioned.

{\bf Acknowledgments:} The authors would like to thank Caio Negreiros, Paolo Piccione and Renato Bettiol for very helpful conversations and comments. 

\section{Variational Setup For the Yamabe Problem}\label{Section1}

Before we present the notions of bifurcation and local rigidity instants, we will describe the set ${\mathcal{R}}^k(M)$ of all {\it $C^k$ Riemannian metrics} on $M$, $m=\dim M\geq 3$, and some properties of the Hilbert-Einstein functional. For further references and details of the concepts and properties presented here see, e.g, \cite{bettiol}, \cite{piccioneproduct} and \cite{schoen}.
	
Let $g_R$ be a fixed auxiliary Riemannian metric on $M$; $g_R$ and its Levi-Civita connection $\nabla^R$ induce naturally Riemannian metrics and connections on all tensor bundles over $M$, respectively. For each $k\geq0$, denote by $\Gamma^k(TM^{\ast}\vee TM^{\ast})$ the space of symmetric $C^k$ sections of $TM^{\ast}\otimes TM^{\ast}$, i.e, symmetric $(0,2)$-tensors of class $C^k$ on $M$. This becomes a Banach space when equipped with the $C^k$ norm $$\left\|\tau\right\|_{C^k}=\displaystyle\max_{j=1,\ldots,k}\left(\displaystyle\max_{x\in M}\left\|(\nabla^R)^j\tau(x)\right\|_R\right),$$ where $\left\|\cdot\right\|_R$ denote the norm induced by $g_R$ on each appropriate space. 
	
Note that the set ${\mathcal{R}}^k(M)$ of all {\it $C^k$ Riemannian metrics} on $M$ is a open convex cone inside $(\Gamma^k(TM^{\ast}\vee TM^{\ast}),\left\|\cdot\right\|_{C^k})$ and, therefore, contractible. Thus, ${\mathcal{R}}^k(M)$ inheriting a natural differential structure. We remark also that ${\mathcal{R}}^k(M)$ is a open set of a Banach space and, thus, we can identify its tangent space with the vector space $\Gamma^k(TM^{\ast}\vee TM^{\ast})$. From now on, we assume that $k\geq3.$
		
For each $g\in{\mathcal{R}}^k(M)$, denote by $\textrm{vol}_g$ the volume form on $M$ (we assume that $M$ is orientable); in this case, $L^2(M,\textrm{vol}_g)$ will denote the usual Hilbert space of real square integrable functions on $M$. Consider over ${\mathcal{R}}^k(M)$ the maps $$\textrm{scal}:{\mathcal{R}}^k(M)\longrightarrow C^{k-2}(M) \ \ \textrm{and} \ \ \textrm{Vol}:{\mathcal{R}}^k(M)\longrightarrow\R,$$ the {\it scalar curvature} and the {\it volume}, that for each Riemannian metric $g\in{\mathcal{R}}^k(M)$ associate, respectively, its scalar curvature $\textrm{scal}(g):M\longrightarrow\R$ and its volume $\textrm{Vol}(g)=\int_M{\textrm{vol}_g}$. Define the {\it Hilbert-Einstein functional} as the function $A:{\mathcal{R}}^k(M)\longrightarrow\R$ given by 
\begin{equation}
\mathcal{A}(g)=\dfrac{1}{\textrm{Vol}(g)}\int_M\textrm{scal}(g)\textrm{vol}_g. \label{HE}
\end{equation}
In order to enunciate some properties of the Hilbert-Einstein functional, we will establish the more appropriate regularity for our manifold of metrics and maps. Denote by ${\mathcal{R}}_{1}^k(M)=\textrm{Vol}^{-1}(1)$ the subset of ${\mathcal{R}}^k(M)$ consisting of metrics of volume one. Note that ${\mathcal{R}}_{1}^k(M)$ is a smooth embedded codimension 1 submanifold of ${\mathcal{R}}^k(M)$. 

%the conformal class of $g$ is %
For each $g\in{\mathcal{R}}^k(M)$, we denote by $[g]\subset {\mathcal{R}}^k(M)$ the set of metrics conformal to $g$. The space $\Gamma^k(TM^{\ast}\vee TM^{\ast})$ defined above induces a differential structure on each conformal class. 

%We wil Later, some Fredholm's conditions must be satisfied. For this, it is enough we introduce the notion of {\it $C^{k,\alpha}$ conformal class} and, therefore, define on the conformal class $g\in{\mathcal{R}}^k(M)$ a Hölder $C^{k,\alpha}(M)$ regularity. 

Let $$[g]_{k,\alpha}=\left\{\phi g;\phi\in C^{k,\alpha}(M), \ \ \phi>0\right\}$$ be the {\it $C^{k,\alpha}(M)$ conformal class} of $g$, which can be identified with the open subset of $C^{k,\alpha}(M)$ formed by the positive functions, which allows $[g]_{k,\alpha}$ to have a natural differential structure. Working with this regularity, the necessary Fredholm's conditions for the second variation of the Hilbert-Einstein functional are satisfied. The set $${\mathcal{R}}_{1}^{k,\alpha}(M,g)={\mathcal{R}}_{1}^k(M)\cap[g]_{k,\alpha}$$is a smooth embedded codimension 1 submanifold  of $[g]_{k,\alpha}$. %Proposition \ref{hilbert} contains well-known facts about $\mathcal{A}$ and its critical points, see, e.g. \cite{piccioneproduct}.

\begin{proposition} [\cite{piccioneproduct}] \label{hilbert}The restriction of the Hilbert-Einstein functional $\mathcal{A}$ to $\mathcal{R}_1^{k,\alpha}(M,g)$ is smooth, and its critical points are constant scalar curvature metrics in $\mathcal{R}_1^{k,\alpha}(M,g)$. 
\end{proposition}

According \cite{piccioneproduct} we also have that given a critical point $g_0\in \mathcal{R}_1^{k,\alpha}(M,g)$ of $\mathcal{A}$, the second variation of $\mathcal{A}$ at $g_0$ can be identified with the quadratic form $$\text{d}^2\mathcal{A}(g_0)(\psi,\psi)=\dfrac{m-2}{2}\int_M((m-1)\Delta_{g_0}\psi-\text{scal}(g_0)\psi)\psi\text{vol}_{g_0},$$ defined on the tangent space at $g_0=\phi g$, given by $$T_{g_0}{\mathcal{R}}_{1}^{k,\alpha}(M,g)\cong\left\{\psi\in C^{k,\alpha}(M);\int_M\frac{\psi}{\phi}\text{vol}_{g_0}=0\right\}.$$

\begin{remark} \emph{We are considering $\Delta_{g}=-\text{div}_g\circ\text{grad}_g$, the Laplacian operator of $(M,g)$, acting on $C^{\infty}(M)$. We observe that, given $\lambda\in\R^+$, one has $\Delta_{\lambda g}=\frac{1}{\lambda}\Delta_g$ and $\text{scal}(\lambda g)=\frac{1}{\lambda}\text{scal}(g)$. Therefore, the negative eigenvalues (respect. positive) of $\Delta_g$ remain negative (respect. positive), that is, we can normalize metrics in order to have unit volume, without changing the spectral theory of the operator $\Delta_g-\dfrac{\text{scal}(g)}{m-1}\cdot\text{Id}$, $m=\dim M$.}
\end{remark} 
 
Now, we will introduce the concepts of {\it bifurcation} and {\it local rigidity} for a 1-parameter family of solutions to the Yamabe problem. Let $$[a,b] \ni t \mapsto g_t\in {\mathcal{R}}^k(M), \ \ k\geq 3,$$ or simply $g_t$, denote a smooth 1-parameter family (smooth path) of Riemannian metrics on $M$, such that each $g_t$ has constant scalar curvature.
 
\begin{definition} An instant $t_{\ast}\in [a,b]$ is a \emph{bifurcation instant} for the family $g_t$ if there exists a sequence $\{t_q\}$ in $[a,b]$ that converges to $t_{\ast}$ and a sequence $\{g_q\}$ in ${\mathcal{R}}^k(M)$ of Riemannian metrics that converges to $g_{t_{\ast}}$ satisfying

\begin{itemize}
\item[(i)] $g_{t_q}\in [g_q]$, but $g_q\neq g_{t_q}$;
\item[(ii)] $\emph{Vol}(g_q)=\emph{Vol}(g_{t_q})$;
\item[(iii)]$\emph{scal}(g_q)$ is constant.
\end{itemize}
\end{definition}

If $t_{\ast}\in [a,b]$ is not a bifurcation instant, it is said that the family $g_t$ is {\it locally rigid} at $t_{\ast}$. More precisely, the family is locally rigid at $t_{\ast}\in [a,b]$ if there exists an open set $U\subset{\mathcal{R}}^k(M)$ containing $g_{t_{\ast}}$ such that if $g\in U$ is another metric with constant scalar curvature and there exists $t\in [a,b]$ with $g_t\in U$ and

\begin{itemize}
\item[\textrm{(a)}] $g\in [g_t]$;
\item[\textrm{(b)}] $\textrm{Vol}(g)=\textrm{Vol}(g_t)$,
\end{itemize} then $g=g_t$.

%In order to identify the occurrence of the above two situations, we must define {\it degeneracy instant} and Morse index of $g_t.$

\begin{definition} An instant $t_{\ast}\in [a,b]$ is a \emph{degeneracy instant} for the family $g_t$ if $\emph{scal}(g_{t_{\ast}})\neq 0$ and $\frac{\emph{scal}(g_{t_{\ast}})}{m-1}$ is an eigenvalue of the Laplacian operator $\Delta_{t_{\ast}}$ of $g_{t_{\ast}}$. 
\end{definition}

\begin{remark} In face of the second variation expression of the Hilbert-Einstein functional the above definition is equivalent to the fact that $g_{t_{\ast}}$ being a degenerate critical point (in the usual sense of Morse theory) of the Hilbert-Einstein functional on ${\mathcal{R}}_{1}^{k,\alpha}(M,g)$.
\end{remark}

The {\it Morse index} $N(g_t)$ of $g_t$ is given by the number of positive eigenvalues of $\Delta_t$ counted with multiplicity less than $\frac{\text{scal}(g_{t_{\ast}})}{m-1}$. For each $t>0$, $N(g_t)$ is a non-negative integer number. In other words, $N(g_t)$ counts the number of directions which the functional decreases, since the second variation is negative definite.

We will now state a general criterion for classifying local rigidity instants for a 1-parameter family $g_t$ of solutions to the Yamabe problem.

\begin{proposition}[\emph{\cite{bettiol}}]\label{criteriorig} Let $g_t$ be a smooth path of metrics of class $C^k$, $k\geq 3$, such that $\emph{scal}(g_t)$ is constant for all $t\in [a,b]$, and let $\Delta_t$ be the Laplacian operator of $g_t$. If $t_{\ast}$ is not a degeneracy instant of $g_t$, then $g_t$ is locally rigid at $t_{\ast}$.
\end{proposition}

\begin{corollary}[\emph{\cite{bettiol}}] \label{criterioauto}Suppose that in addition to the hypotheses of Proposition \ref{criteriorig}, there exists an instant $t_{\ast}$ when $\frac{\emph{scal}(g_{t_{\ast}})}{m-1}$ is less than the first positive eigenvalue $\lambda_1(t_{\ast})$ of $\Delta_{t_{\ast}}$. Then $g_{t_{\ast}}$ is a local minimum for the Hilbert-Einstein functional in its conformal class. In particular, $g_t$ locally rigid at $t_{\ast}.$
\end{corollary}

\begin{remark}\emph{Note that, in fact, the Corollary \ref{criterioauto} is an immediate consequence from Proposition \ref{criteriorig}, since Morse index of the non-degenerate critical point $g_{t_{\ast}}$ is $N(g_{t_{\ast}})=0.$ In addition, we remark that, if $t_{\ast}$ is a bifurcation instant for $g_t$, then $t_{\ast}$ is necessarily a degeneracy instant for $g_t$. However, the reciprocal is not true in general.} 
\end{remark}

The next result provides a sufficient condition to determinate if a degeneracy instant is a bifurcation instant, given in terms of a jump in the Morse index when passing a degeneracy instant. 

\begin{proposition}\emph{\cite{bettiol}} \label{criterioindice}Let $g_t$ a smooth path of metrics of class $C^k$, $k\geq 3$, such that $\emph{scal}(g_t)$ is constant for all $t\in [a,b]$, $\Delta_t$ the Laplacian $g_t$ and $N(g_t)$ the Morse index of $g_t$. Assume that $a$ and $b$ are not degeneracy instants for $g_t$ and $N(g_a)\neq N(g_b)$. Then, there exists a bifurcation instant $t_{\ast}\in\, ]a,b[$ for the family $g_t$.
\end{proposition}

\section{Generalized Flag Manifolds}\label{flag}
	
%	Let $G$ be a compact semissimple Lie group and let $x\in \lie$ be a regular element, i.e, its centralizer is equal to a maximal toral subalgebra of $\lie=\text{Lie}(G)$, Lie algebra of $G$. The {\it adjoint orbit} of $x$ is the set $M=G\cdot x=\Ad(G)x \subset \lie$. Let $$K=G_x=\{g\in G : \Ad(g)x=x\}\subset G$$be isotrpy subgroup of $x$. Then $M$ is diffeomorphic to the left coset $G/K$. 
%	
%	Let $x_0$ be the point which corresponds to the class $1\cdot K\in G/K$ ($1\in G$ is the identity element of $G$). If we take $T_{x_0} = \overline{\exp \R x_0}$, it is well known that $T_{x_0}$ is a toral subgroup of $G$ (diffeomorphic to the Lie group $S^1 \times S^1\times\ldots\times S^1$). Furthermore, one has $K=K_x=C(T_x)$, where $C(T_x)$ is the centralizer of $T_x$ in $G$. If $T_x$ is maximal, $C(T_x)=T_x$. 
%	
%Using the above notations, we can introduce the concept of {\it generalized flag manifold} which we will consider in this work.	
	
\begin{definition}	
Let $G$ be a compact semissimple Lie group with Lie algebra $\lie$. A \emph{generalized flag manifold} is the adjoint orbit of a regular element in the Lie algebra $\lie$. Equivalently, a generalized flag manifold is a homogeneous space of the form $G/C(T)$, where $T$ is a torus in $G$. When $T$ is maximal, it is said that $G/C(T)=G/T$ is a \emph{maximal flag manifold}. 
\end{definition}

\begin{example} The maximal flag manifolds of a classical Lie group are given by:
\begin{enumerate}
\item[{\bf A:}] $SU(n+1)/S(U(1)\times\ldots\times U(1))=SU(n+1)/S(U(1)^{n+1})$; 
\item[{\bf B:}] $SO(2n+1)/U(1)\times\ldots\times U(1)=SO(2n+1)/U(1)^n$;
\item[{\bf C:}] $Sp(n)/U(1)\times\ldots\times U(1)=Sp(n)/U(1)^n$;
\item[{\bf D:}] $SO(2n)/U(1)\times\ldots\times U(1)=SO(2n)/U(1)^n$;
\item[{\bf G:}] $G_2/U(1)\times U(1)$;
\end{enumerate}

\end{example}

%The Arvanitoyeorgos' work \cite{Arvanit} also introduce a complete classification of partial flag manifolds.

%Thereafter, we will define generalized flag manifold from a complex Lie algebra.
Let us give a Lie theoretical description of generalized flag manifolds. Let $\lie^{\mathbb{C}}$ be a complex semissimple Lie algebra. Given a Cartan subalgebra $\mathfrak{h}^{\mathbb{C}}$ of $\lie^{\mathbb{C}}$,
denote by $R$ the set of $roots$ of $\lie^{\mathbb{C}}$ with relation to $\mathfrak{h}^{\mathbb{C}}$. Consider the decomposition $$\lie^{\mathbb{C}}=\mathfrak{h}^{\mathbb{C}}\oplus\sum_{\alpha\in R}\lie_{\alpha}^{\mathbb{C}},$$ where $\lie_{\alpha}^{\mathbb{C}}=\{X\in\lie^{\mathbb{C}};\forall H\in \mathfrak{h}^{\mathbb{C}}, \left[X,H\right]=\alpha(H)X\}$. Let $R^+\subset R$ be a choice of positive roots, $\Sigma$ the correspondent system of simple roots and $\Theta$ a subset of $\Sigma$. Will be denoted by $\left\langle \Theta\right\rangle$ the span of $\Theta$, $R_M=R\setminus \left\langle \Theta\right\rangle$ the set of the {\it complementary roots} and by $R_{M^+}$ the set of positive complementary roots.

A Lie subalgebra $\mathfrak{p}$ of $\lie^{\mathbb{C}}$ is called {\it parabolic} if it contains a Borel subalgebra of $\lie^{\mathbb{C}}$ (i.e., a maximal solvable Lie subalgebra of $\lie^{\mathbb{C}}$). Take $$\mathfrak{p}_{\Theta}=\mathfrak{h}^{\mathbb{C}}\oplus\sum_{\alpha\in \left\langle \Theta\right\rangle^+}\lie_{\alpha}^{\mathbb{C}}\oplus\sum_{\alpha\in \left\langle \Theta\right\rangle^+}\lie_{-\alpha}^{\mathbb{C}}\oplus\sum_{\beta\in \Pi_{M}^+}\lie_{\beta}^{\mathbb{C}}$$ the canonical parabolic subalgebra of $\lie^{\mathbb{C}}$ determined by $\Theta$ which contains the Borel subalgebra $\mathfrak{b}=\mathfrak{h}^{\mathbb{C}}\oplus\displaystyle\sum_{\beta\in \Pi_+}\lie_{\beta}^{\mathbb{C}}$.

The generalized flag manifold $\mathbb{F}_{\Theta}$ associated with $\lie^{\mathbb{C}}$ is defined as the homogeneous space $$\mathbb{F}_{\Theta}=G^{\mathbb{C}}/P_{\Theta},$$where $G^{\mathbb{C}}$ is the complex compact connected simple Lie group with Lie algebra $\lie$ and $P_{\Theta}=\{g\in G^{\mathbb{C}}; \Ad(g)\mathfrak{p}_{\Theta}=\mathfrak{p}_{\Theta}\}$ is the normalizer of $\mathfrak{p}_{\Theta}$ in $G^{\mathbb{C}}$. 

Let $G$ be the real compact form of $G^{\mathbb{C}}$ corresponding to $\lie$, i.e, $G$ is the connected Lie Group with Lie algebra $\lie$, the compact real form of $\lie^{\mathbb{C}}$. The subgroup $K_{\Theta}=G\cap P_{\Theta}$ is the centralizer of a torus. Furthermore, $G$ acts transitively on $\mathbb{F}_{\Theta}$. Since $G$ is compact, we have that $\mathbb{F}_{\Theta}$ is a compact homogeneous space, that is, $$\mathbb{F}_{\Theta}=G^{\mathbb{C}}/P_{\Theta}=G/G\cap P_{\Theta}=G/K_{\Theta},$$ accordingly characterization given above from adjoint orbits of a regular element of $\lie$. 

There are two classes of generalized flag manifolds. The first occurs when $\Theta=\emptyset$. Therefore, the parabolic subalgebra is given by  $\mathfrak{p}_{\Theta}=\mathfrak{h}^{\mathbb{C}}\oplus\displaystyle\sum_{\beta\in \Pi_+}\lie_{\beta}^{\mathbb{C}}$, that is, is equal to the Borel subalgebra of $\lie^{\mathbb{C}}$ and $T=P_{\Theta}\cap G$ is a maximal torus. In this case, $\mathbb{F}_{\Theta}=G^{\mathbb{C}}/P_{\Theta}=G/G\cap P_{\Theta}=G/T$ is called {\it maximal} flag manifold. When $\Theta\neq\emptyset$, $\mathbb{F}_{\Theta}$ is called {\it generalized} flag manifold. 

We now consider a {\it Weyl basis} for $\lie^{\mathbb{C}}$ given by $\{X_{\alpha}\}_{\alpha\in R}\cup\{H_{\alpha}\}_{\alpha\in \Sigma}\subset \lie^{\mathbb{C}}, \ \ X_{\alpha}\in \mathfrak{g}_{\alpha}^{\mathbb{C}}$. From this basis we determine a basis for $\lie$, the compact real form of $\lie^{\mathbb{C}}$, putting $$\lie=\textrm{span}_{\R}\{\sqrt{-1}H_{\alpha}, A_{\alpha}, S_{\alpha}\},$$ with $A_{\alpha}=X_{\alpha}-X_{-\alpha}$, $S_{\alpha}=\sqrt{-1}(X_{\alpha}+X_{-\alpha})$ ($A_{\alpha}=S_{\alpha}=0$ if $\alpha\notin R$) and $H_{\alpha}\in \mathfrak{h}^{\mathbb{C}}$ is such that $\alpha(\cdot)=\left\langle H_{\alpha},\cdot\right\rangle$, $\alpha\in R$. The {\it structure constants} of this basis are determined by the following relations

$$\left\{
  \begin{array}{rcl} 
  \left[A_{\alpha},S_{\beta}\right]&=& m_{\alpha,\beta}A_{\alpha+\beta}+m_{-\alpha,\beta}A_{\alpha-\beta}\\
  \left[S_{\alpha},S_{\beta}\right]&=& -m_{\alpha,\beta}A_{\alpha+\beta}-m_{\alpha,-\beta}A_{\alpha-\beta}\\
  \left[A_{\alpha},S_{\beta}\right]&=& m_{\alpha,\beta}S_{\alpha+\beta}+m_{\alpha,-\beta}S_{\alpha-\beta}\\
  \end{array}
  \right.,$$
	
	$$\left\{
  \begin{array}{rcl}  
  \left[\sqrt{-1}H_{\alpha},A_{\beta}\right]&=& \beta(H_{\alpha})S_{\beta}\\
 \left[\sqrt{-1}H_{\alpha},S_{\beta}\right]&=& -\beta(H_{\alpha})A_{\beta}\\
  \left[A_{\alpha},S_{\alpha}\right]&=& 2\sqrt{-1}H_{\alpha}\\
  \end{array}
  \right.,$$
	where $m_{\alpha,\beta}$ is such that $\left[X_{\alpha},X_{\beta}\right]=m_{\alpha,\beta}X_{\alpha+\beta}$, with $m_{\alpha,\beta}=0$ if $\alpha+\beta\notin R$ and $m_{\alpha,\beta}=-m_{-\alpha,-\beta}$. We remark that this basis is $-B$-orthogonal and $-B(A_{\alpha},A_{\alpha})=-B(S_{\alpha},S_{\alpha})=2$, where $B$ Cartan-Killing form of $\lie^{\mathbb{C}}$ (the Cartan-Killing form of $\lie^{\mathbb{C}}$ restricted to $\lie$ coincides with the Cartan-Killing form of its compact form $\lie$ and $\h_{\R}=\textrm{span}_{\R}\{\sqrt{-1}H_{\alpha}\}_{\alpha\in R}$ is a Cartan subalgebra of $\lie$). Moreover, if $\mathfrak{q}$ and $\mathfrak{s}$ are the subspaces spanned by $\{H_{\alpha},A_{\alpha}\}$ and $\{S_{\alpha}\}$ respectively, one has $$\left[\mathfrak{q},\mathfrak{q}\right]\subset\mathfrak{q}, \ \  \left[\mathfrak{q},\mathfrak{s}\right]\subset\mathfrak{s}, \ \ \left[\mathfrak{s},\mathfrak{s}\right]\subset\mathfrak{q}.$$
	
The above construction of the Weyl basis can be found, e.g., in \cite{SM}, p. 334. This will be the basis that we use when we dealing with flag manifolds in this work.

\subsection{Isotropy Representation}  

From the notion of {\it isotropy representation} we can determine invariant metrics on certain homogeneous spaces as follows.
 
Let $G\times M\longrightarrow M$ be a defferentiable and transitive action of a Lie group $G$ on the homogeneous space $(M,m)$ endowed with a $G$-invariant metric $m$. Given $x\in G$, let $K=G_x$ be the isotropy subgroup of $x$. The {\it isotropy representation} of $K$ is the homomorphism $g\in K\mapsto dg_x \in \textrm{Gl}(T_xM)$. Note that $m_{x_0}$ is a inner product on $T_{x_0}(G/K)$, invariant by such representation, where $x_0=1\cdot K$.

A homogeneous space $G/K$ is {\it reducible} if $G$ has a Lie algebra $\lie$ such that $$\lie=\mathfrak{k}\oplus\mathfrak{m},$$ with $\Ad(K)\mathfrak{m}\subset \mathfrak{m}$ and $\mathfrak{k}$ Lie algebra of $K$. If $K$ is compact, this decomposition always exists, namely, if we take $\mathfrak{m}=\mathfrak{k}^{\bot}$, $-B$-orthogonal complement to $\mathfrak{k}$ in $\lie$, where $B$ is the Cartan-Killing form of $\lie$.

If $G/K$ is reducible the isotropy representation of $K$ is equivalent to $\Ad |_K$, the restriction of the adjoint representation of $G$ to $K$, that is,  $$j(k)=\Ad(k)|_{\mathfrak{m}}, \forall k\in K.$$

A representation of a compact Lie group $K$ is always orthogonal (preserves inner product) on the representation space. We can conclude that every reductive homogeneous space $G/K$ has a $G$-invariant  metric, since such a metric is completely determined by an inner product on the tangent space at the origin $T_{x_0}(G/K)$.  

We remark that the set of all $G$-invariant metrics on $G/K$ are in 1-1 correspondence the set of inner products $\left\langle, \right\rangle$ on $\mathfrak{m}$, invariant by $\Ad(k)$ on $\mathfrak{m}$, for each $k\in K$, that is, $$\left\langle \Ad(k)X, \Ad(k)Y \right\rangle=\left\langle X, Y \right\rangle, \forall X,Y\in\mathfrak{m}, k\in K.$$

In the case of generalized flag manifolds the isotropy representation of $K$ leaves $\m$ invariant, i.e, $\Ad(K)\mathfrak{m}\subset \mathfrak{m}$ and decomposes it into irreducible submodules $$\m=\m_1\oplus\m_2\ldots\oplus\m_n,$$ and these submodules are inequivalent to each other. The submodules $\m_i$ are called {\it isotropy summands}. 

It follows that a $G$-invariant metric $g$ on $G/K$ is represented by a inner product $$g_{1\cdot K}=t_1Q|_{\m_1}+t_2Q|_{\m_2}+\ldots +t_nQ|_{\m_n}$$on $\m$, with $t_i$ positive constants and $Q$ is (the extension of) a inner product on $\m$, $\Ad(K)$-invariant. 

In particular, if $Q=(-B)$, with $B$ Cartan-Killing form of $G$ and $t_i=1$ for all $1\leq i\leq n$ above, the $G$-invariant metric $g$ on $G/K$ represented by the inner product $$g_{1\cdot K}=(-B)|_{\m_1}+(-B)|_{\m_2}+\ldots +(-B)|_{\m_n}$$on $\m$ is called {\it normal metric}.

\begin{remark} Let $M=G/T$ be a maximal flag manifold of a compact simple Lie group $G$. Then, the isotropy representation of $M$ decomposes into a discrete sum of $2$-dimensional pairwise non-equivalent irreducible $T$-submodules $\m_{\alpha}$ as follows: $$\m=\sum_{\alpha\in R^+}\m_{\alpha}.$$The number of these submodules is equal to the cardinality of $R^+$, the set of positive roots of $\mathfrak{g}$.
\end{remark} 

\begin{table}[htb]
\caption{The number of isotropy summands for maximal flag manifolds $G/T$}
\centering
\begin{tabular}{lcl}
\hline
Maximal flag manifold $G/T$ & Number of roots $|R|$ &$\m=\oplus_{j=1}^l$\\ 
\hline 
\vspace{0.3cm}
$SU(n+1)/T^{n},n\geq 1$&$n(n+1)$ & $l=n(n+1)/2$\\[1pt] \vspace{0.3cm}
$SO(2n+1)/T^n,n\geq 2$&$2n^2$ & $l=n^2$\\ \vspace{0.3cm} 
$Sp(n)/T^n,n\geq 3$&$2n^2$ & $l=n^2$\\ \vspace{0.3cm}
$SO(2n)/T^n,n\geq 4$&$2n(n-1)$ & $l=n(n-1)$\\ \vspace{0.3cm} 
$G_2/T$&$12$ & $l=6$\\ \vspace{0.3cm}
\end{tabular}
\end{table}
% \newpage
%\begin{remark}
%For $n=1$ the maximal flag $SU(n+1)/T^{n}$ is $SU(2)/S(U(1)\times U(1))\cong \mathbb{CP}^1,$ which is an isotropy irreducible Hermitian symmetric space.
%\end{remark}

\section{Laplacian and Canonical Variation of Riemannian Submersions With Totally Geodesic Fibers}
\label{cvariation}

%% parei aqui

The Lie theoretic description of flag manifolds is applied in order to obtain the horizontal and vertical distributions for each homogeneous fibration whose the total space represents a class of maximal flag manifold $G/T$ provided with a {\it normal} metric $g$, where $G$ is a compact simple Lie group and $T\subset G$ is a maximal torus in $G$.

We denote by $\Delta_{g}=-\textrm{div}_g\circ\textrm{grad}_g$ the Laplacian operator of $(M,g)$ acting on $C^{\infty}(M)$. The operator $\Delta_{g}$, densely defined on $L^2(M,\textrm{vol}_g)$, is symmetric (hence closable), non-negative has compact resolvent. Furthermore, it is well-known that $\Delta_g$ is essentially self-adjoint with this domain. We denote its unique self-adjoint extension also by $\Delta_g$. Analogously, let $\Delta_k$ be the unique self-adjoint extension of the Laplacian of the fiber $(F,k)$, where $k$ is the metric induced by $(M,g)$ on $F$.

\begin{definition} [\cite{bergery}] The \emph{vertical Laplacian} $\Delta_v$ acting on $L^2(M,\emph{vol}_g)$ is the operator defined at $p\in M$ by $$(\Delta_v\psi)(p)=(\Delta_k\psi|_{F_p})(p),$$ and the \emph{horizontal Laplacian} $\Delta_h$, acting on the same space, is defined by the difference $$\Delta_h=\Delta_g-\Delta_v.$$ 
\end{definition}

Both $\Delta_h$ and $\Delta_v$ are non-negative self-adjoint unbounded operators on $L^2(M,\textrm{vol}_g)$, but in general, are not elliptic (unless $\pi$ is a covering). We now consider the spectrum of such operators. 

As remarked above, $\Delta_g$ is non-negative and has compact resolvent, that is, its spectrum is non-negative and discrete. Since the fibers are isometrics, $\Delta_v$ also has discrete spectrum equal to the fibers. However, the spectrum of $\Delta_h$ need not be discrete.

%\begin{theorem} [\cite{bergery}]\label{bergery} If the fibers of the Riemannian submersion $\pi:M\longrightarrow B$ are totally geodesics the operators $\Delta_g$, $\Delta_v$ e $\Delta_{g'}$ commute with each other.
%\end{theorem}

%When $\Delta_g$ and $\Delta_v$ commute, we have the following decomposition of $L^2(M,g)$.

%\begin{theorem} [\cite{bergery}] The Hilbert space $L^2(M,g)$ admits a Hilbert basis consisting of simultaneous eigenfunctions of $\Delta_g$ e $\Delta_v$.
%\end{theorem}
\subsection{Canonical Variation}\label{Section4}

%If we fix Riemannian submersion with totally geodesic fibers, it is possible to define a 1-parameter family $g_t$, $t>0$, of other such Riemannian submersions by scaling the original metric of the total space in the direction of the fibers. 

Let us recall the definition of {\it the canonical variation}. Its will be fundamental in the determination of bifurcation and local rigid instants in our work.

\begin{definition}[\cite{bergery}]\label{variation} Let $F\cdots (M,g)\stackrel{\pi}{\rightarrow} B$ be a Riemannian submersion with totally geodesic fibers. Consider the 1-parameter family of Riemannian submersions given by \linebreak $\left\{F\cdots (M,g_t)\stackrel{\pi}{\rightarrow} B, t>0\right\}$, where $g_t\in {\mathcal{R}}^k(M)$ is defined by 
$$g_t(v,w)=\left\{
  \begin{array}{rcl} 
  t^2g(v,w),&\mbox{se} &v,w \ \ \mbox{are verticals}\\
  0,&\mbox{se} &v \ \ \mbox{is vertical and} \ \ w \ \ \mbox{is horizontal}\\
  g(v,w),&\mbox{se} &v,w \ \ \mbox{are horizontals}.\\
  \end{array}
  \right.$$

Such family of Riemannian submersions is called the \emph{canonical variation} of \linebreak$F\cdots (M,g)\stackrel{\pi}{\rightarrow} B$ or, for simplicity, we may also refer to the family of \emph{total spaces} of these submersions, i.e., the Riemannian manifolds $(M,g_t)$, as the canonical variation of $(M,g)$.
	\end{definition}
	
\begin{proposition} [\cite{bergery}] The family $\left\{F\cdots (M,g_t)\stackrel{\pi}{\rightarrow} B, t>0\right\}$ of Riemannian submersions has totally geodesic fibers, for each $t>0$. Furthermore, its fibers are isometrics to $(F,t^2k)$, where $(F,k)$ is the original fiber of $\pi:M\longrightarrow B$. 
\end{proposition}

\begin{remark}
\emph{Note that, for $a\neq b$, the metrics $g_a$ e $g_b$ are not conformal. Furthermore, for each $t>0$, $g_t$ is the unique Riemannian metric that satisfy the conditions of Definition \ref{variation}.}

\end{remark}

The following result shows how to decompose $\Delta_t$ in terms of the vertical and horizontal Laplacians.

\begin{proposition} [\cite{bettiol}] Let $\Delta_t$ the Laplacian of $(M,g_t)$. Then 

\begin{equation}
\Delta_t=\Delta_h+\frac{1}{t^2}\Delta_v=\Delta_g+(\frac{1}{t^2}-1)\Delta_v. \label{deltat-eq}
\end{equation}
\end{proposition}

\begin{corollary} [\cite{bettiol}] \label{deltat}For each $t>0$, the following inclusion holds 
\begin{equation}
\sigma(\Delta_t)\subset\sigma(\Delta_g)+(\frac{1}{t^2}-1)\sigma(\Delta_v), \label{inclusao}
\end{equation}
where $\sigma(\Delta_t)$, $\sigma(\Delta_g)$ and $\sigma(\Delta_v)$ are the respective spectrum of $\Delta_t$, $\Delta_g$ e $\Delta_v.$ Since the above spectra are discrete, this means that every eigenvalue $\lambda(t)$ of $\Delta_t$ is of the form
\begin{equation}
\lambda^{k,j}(t)=\mu_k+(\frac{1}{t^2}-1)\phi_j, \label{kj}
\end{equation}
for some $\mu_k\in\sigma(\Delta_g)$ and $\phi_j\in\sigma(\Delta_v).$ 
\end{corollary}
%{\bf Proof:} By Theorem \ref{bergery}, $\Delta_g$ and $\Delta_v$ commute. Hence, by the Spectral Theorem, such operators are simultaneously diagonalizable, in the sense that there exists a unitary operator $U$ of $L^2(M,\textrm{vol}_g)$ such that $U\Delta_gU^{-1}=f_M$ and $U\Delta_vU^{-1}=f_v$ are multiplication operators by functions $f_M$ and $f_v$ respectively. For such multiplication operators, the spectrum $\sigma(T_f)$ is the essential range, $\textrm{ess.Im}(f)(\subset\overline{\textrm{Im}(f)})$, of $f$. Then, by the expression of $\Delta_t$ in \eqref{deltat},{\small $$\sigma(\Delta_t)=\textrm{ess.Im}(f_M+(\frac{1}{t^2}-1)f_v)\subset\overline{\textrm{ess.Im}(f_M)+\textrm{ess.Im}((\frac{1}{t^2}-1)f_v)}=\overline{\sigma(\Delta_g)+(\frac{1}{t^2}-1)\sigma(\Delta_v)}.$$}Since both spectra are discrete, we may remove the closure and the inclusion \eqref{inclusao} is proved. \cqd   

\begin{corollary} \label{primeiroaut1}If $\lambda_1(t)$ is the first positive eigenvalue of $\Delta_t$, then $$\lambda_1(t)\geq \mu_1, \forall \ \ 0<t\leq 1,$$ where $\mu_1$ is the first positive eigenvalue of $\Delta_g$.
\end{corollary}
{\bf Proof:} Since $\lambda_1(t)=\mu_k+(\frac{1}{t^2}-1)\phi_j$, for some $\mu_k\in \sigma(\Delta_g)$ and $\phi_j\in \sigma(\Delta_v)$, if $0<t\leq 1$, $(\frac{1}{t^2}-1)\phi_j\geq 0$ and $\lambda_1(t)=\mu_k+(\frac{1}{t^2}-1)\phi_j\geq\mu_k\geq\mu_1,$ where $\mu_1\in\sigma(\Delta_g)$ is the first positive eigenvalue of the operator $\Delta_g$. \cqd

We also remark that not all possible combinations of $\mu_k$ and $\phi_j$ in the expression \eqref{kj} give rise to an eigenvalue of $\Delta_t$. In fact, this only happens when the total space $(M,g)$ of the submersion is a Riemannian product. To determine which combinations are allowed in general is difficulty and depends on the global geometry of the submersion. Moreover, the ordering of the eigenvalues of $\Delta_t$ may change with  $t$. 

We have the following useful properties of the spectrum of the Laplacian operator $\Delta_t$ of $g_t$.

\begin{proposition} [\cite{pacificjournal}] Using the same notations above, one has that $$\sigma(\Delta_B)\subset\sigma(\Delta_t),$$ for all $t>0.$
\end{proposition}
%{\bf Proof:} For each $\psi:B\longrightarrow\R$ and its lift $\widetilde{\psi}:=\psi\circ\pi$, 

%\begin{equation}
%\Delta_g\widetilde{\psi}=(\Delta_B\psi)\circ\pi+g(\textrm{grad}_g\widetilde{\psi},\vec{H}), \label{delta}
%\end{equation}
%where $\vec{H}$ is the {\it mean curvature} vector field of the fibers. Since we assumed the fibers of $\pi$ are totally geodesics, $\vec{H}\equiv 0$ over the fibers. It follows that, if $\psi$ is a eigenfunction of $\Delta_B$, then its lift $\widetilde{\psi}$ is an eigenfunction of $\Delta_g$ with the same eigenvalue (and constant along the fibers) and therefore $$\sigma(\Delta_B)\subset\sigma(\Delta_g).$$ Since the fibers of $\pi$ are totally geodesics with respect to $g_t$, the above inclusion above holds when $\Delta_g$ is replaced with $\Delta_t$, i.e, $$\sigma(\Delta_B)\subset\sigma(\Delta_t), \ \ t>0.$$
%\cqd

We remark that the spectrum of $\Delta_h$, $\sigma(\Delta_h)$, contains but not coincides with the spectrum of the basis $B$. In fact, if $\overline{f}$ is a $C^{\infty}$ function on the basis $B$, then $$(\Delta_{g'}\overline{f})\circ\pi=\Delta_g(\overline{f}\circ\pi)=\Delta_h(\overline{f}\circ\pi),$$ where $\Delta_{g'}$ is the Laplacian operator acting on functions in $C^{\infty}(B,g')$.

\begin{corollary} \label{primeiroaut2}Denoting by $\beta_1$ the first positive eigenvalue of $\Delta_B$ and by $\lambda_1(t)$ the first positive eigenvalue of $\Delta_t$, the following inequality holds $$\lambda_1(t)\leq \beta_1, \ \ \forall \, t>0.$$
\end{corollary}
  
As consequence of the Corollaries \ref{primeiroaut1} and \ref{primeiroaut2}, we have that $\mu_1\leq\lambda_1(t)\leq \beta_1$, where $\mu_1$ the first positive eigenvalue of the Laplacian $\Delta_g$ on the total space $M$. 

When $j=0$ in the expression \eqref{kj}, if $\lambda^{k,0}(t)=\mu_k\in\sigma(\Delta_g)$ remains an eigenvalue of $\Delta_t$ for $t\neq 1$, such eigenvalues will be called {\it constant eigenvalue of $\Delta_t$}, since they are independent of $t$. We stress that $\lambda^{k,0}(t)$ is not necessarily a constant eigenvalue of $\Delta_t$ for all $k$. A simple criterion to determinate when $\lambda^{k,0}(t)\in\sigma(\Delta_t)$ is given in the following:

\begin{proposition} [\cite{pacificjournal}] \label{constanteigen}$\mu_k=\lambda^{k,0}(t)\in\sigma(\Delta_t)$ for $t\neq 1$ $\Leftrightarrow$ $\mu_k\in\sigma(\Delta_B)$.
\end{proposition}

%\begin{corollary} \label{constanteig}$\mu_k=\lambda^{k,0}(t)\in\sigma(\Delta_t)$ for $t\neq 1$ $\Leftrightarrow$ $\mu_k\in\sigma(\Delta_B)$.
%\end{corollary}
%{\bf Proof:} Fixed $t\neq 1$, since $$\lambda^{k,j}(t)=\mu_k+(\frac{1}{t^2}-1)\phi_j \ \ \text{and} \ \ \Delta_t=\Delta_g+(\frac{1}{t^2}-1)\Delta_v,$$ we have that $\lambda^{k,j}(t)\in\sigma(\Delta_t)$ if and only if there exists $\widetilde{\psi}\in C^{\infty}(M)$ such that $$\Delta_g\widetilde{\psi}=\mu_k\widetilde{\psi} \ \ \text{and}\ \ \Delta_v\widetilde{\psi}=\phi_j\widetilde{\psi},$$ with $\mu_k\in\sigma(\Delta_g)$ and $\phi_j\in\sigma(\Delta_v)$. In particular, if we have $\lambda^{k,0}(t)=\mu_k+(\frac{1}{t^2}-1)0\in\sigma(\Delta_t)$, there exists $\widetilde{\psi}$ such that $$\Delta_g\widetilde{\psi}=\mu_k\widetilde{\psi} \ \ \text{and} \ \ \Delta_v\widetilde{\psi}=0\cdot\widetilde{\psi}=0.$$ However, $\Delta_v\widetilde{\psi}=0$ implies that $\widetilde{\psi}$ is constant along the fibers. Therefore, by Proposition 7, there exists a eigenfunction $\psi:B\longrightarrow\R$ in $C^{\infty}(B)$ such that $\Delta_B\psi=\mu_k\psi$ and $\widetilde{\psi}=\psi\circ\pi$. It follows that $\mu_k\in\sigma(\Delta_B)$. The reciprocal results from the fact that $\sigma(\Delta_B)\subset\sigma(\Delta_t)$. \cqd 
\subsubsection{Homogeneous Fibration}\label{Section2.2}
In our main result, the canonical variations $g_t$, obtained from a Riemannian submersion with totally geodesic fibers, are {\it homogeneous metrics}, which are trivial solutions to the Yamabe problem, since every homogeneous metric has constant scalar curvature. This makes these metrics good candidates for admitting other solutions in their conformal class. 

The homogeneous fibrations are obtained from the following construction. Let $K\subsetneq H\subsetneq G$ be compact connected Lie groups, such that $\dim K/H\geq 2$. Consider the natural fibration
$$
\begin{array}{cccc} 
\pi :&\! G/K &\!\longrightarrow &\!G/H\\
     &\! \alpha K &\!\mapsto &\! \alpha H,
		\end{array}
		$$
with fibers $H/K$ $(\textrm{em} \ \ eK)$ and structural group $H$. More precisely, $\pi$ is the associated bundle with fiber $H/K$ to the $H$-principal bundle $p:G\longrightarrow G/H$. 	

Let $\lie$ be the Lie algebra of $G$ and $\mathfrak{h}\supset\mathfrak{k}$ the Lie algebras of $H$ and $K$, respectively. Given a inner product on $\lie$, determined by the Cartan-Killing form $B$ of $\lie$, since $K,H$ and $G$ are compacts, we can consider a $\Ad_G(H)$-invariant orthogonal complement $\mathfrak{q}$ to $\mathfrak{h}$ in $\lie$, i.e., $[\h,\mathfrak{q}]\subset\mathfrak{q}$, and a $\Ad_G(K)$-invariant orthogonal complement $\mathfrak{p}$ to $\sub$ in $\h$, i.e., $[\sub,\mathfrak{p}]\subset\mathfrak{p}$. It follows that $\mathfrak{p}\oplus\mathfrak{q}$ is a $\Ad_G(K)$-invariant orthogonal complement to $\sub$ in $\lie$. 

The $\Ad_G(H)$-invariant inner product $(-B)|_{\mathfrak{q}}$ on $\mathfrak{q}$ define a $G$-invariant Riemannian metric $\breve{g}$ on $G/H$, and the inner product $(-B)|_{\mathfrak{p}}$, $\Ad_G(K)$-invariant on $\mathfrak{p}$, define on $H/K$ a $H$-invariant Riemannian metric $\hat{g}$ on $H/K$. Finally, the orthogonal direct sum of these inner products on $\mathfrak{p}\oplus\mathfrak{q}$ define a $G$-invariant metric $g$ on $G/K$, determined by 
\begin{equation}
g(X+V,Y+W)_{eK}=(-B)|_{\mathfrak{q}}(X,Y) + (-B)|_{\mathfrak{p}}(V,W), \label{metric} 
\end{equation}
for all $X,Y\in\mathfrak{q}$ and $V,W\in\mathfrak{p}$; $g$ is a normal homogeneous metric on $G/K$ and the $(-B)$-orthogonal direct sum $\m=\p\oplus\q$ is the isotropy representation of $K$.

\begin{theorem}[(\cite{besse}, p. 257)] The map $\pi:(G/K,g)\longrightarrow (G/H,\breve{g})$ is a Riemannian submersion with totally geodesic fibers and isometric to $(H/K,\hat{g})$.
\end{theorem} 

If we take for each $t>0$ the metric $g_t$ that corresponds to the inner product on $\mathfrak{p}\oplus\mathfrak{q}$
\begin{equation}
 \left\langle \cdot,\cdot\right\rangle_t=(-B)|_{\mathfrak{q}}+(-t^2B)|_{\mathfrak{p}}, \label{gt} 
\end{equation}
where $B$ is the Cartan-Killing form of $\lie$, we have that the map $\pi_t:(G/K,g_t)\longrightarrow (G/H,\breve{g})$ is a Riemannian submersion with totally geodesic fibers and isometric to $H/K$ provided with the induced metric $-t^2B|_{\mathfrak{p}}$. Hence, we obtain the canonical variation of the original homogeneous fibration $\pi:(G/K,g)\longrightarrow (G/H,\breve{g})$.

\subsection{Canonical Variations of Normal Metrics on Maximal Flag Manifolds}\label{Section2.3}

In this section we will present the construction of the 1-parameter families of homogeneous metrics from normal metric on maximal flag manifolds.

\begin{remark} We must consider the following:

\begin{flushleft}

1) From now on, the total spaces of the original homogeneous fibrations are endowed with a normal homogeneous metric induced by the Cartan-Killing form of the associated complex classical simple Lie algebra. 
\end{flushleft}
%\item[2] \emph{For each complex classical simple Lie algebra, will be used its real compact form with a Weyl basis in order to determine its respective flag manifold;}
\begin{flushleft}
2) Given $K\subsetneq H\subsetneq G$, compact connected Lie groups, with Lie algebras $\sub,\h$ and $\lie$ respectively, the homogeneous fibrations $\pi:G/K\longrightarrow G/H$ are such that the basis space $G/H$ is a Hermitian symmetric space with irredudcible isotropy representation and a metric induced by the restriction of the Cartan-Killing form of $\lie$ to its isotropy representation (horizontal distribution) and the fiber $H/K$ (with $\dim H/K\geq 2$) is some maximal flag manifold, i.e, $K$ is a maximal torus contained in $H$, and $H/K$ is provided with the normal metric induced by the Cartan-Killing form of $\lie$ restricted to its isotropy representation (vertical distribution).
\end{flushleft}

\end{remark}

\begin{proposition} [\cite{SM}]Let $\lie$ and $\h$ be complex classic simple Lie algebras such that $\h\subset\lie$, and let $\mathfrak{t}$ be a Cartan subalgebra for $\lie$ with Cartan-Killing form $B$. If $R$ is the root system of $\mathfrak{t}$ associated with $\lie$, then there exists a subset $R'$ of $R$ such that
\begin{itemize}
\item[\emph{(a)}] The subspace $\mathfrak{t}(R')$ of $\h$ spanned by the dual vectors $H_{\alpha}, \alpha\in R'$, is a Cartan subalgebra of $\h$ and $R'$ is a system of roots of $\mathfrak{t}(R')$.
\item[\emph{(b)}] The Lie algebra $\h$ has a $(-B)$-orthogonal decomposition into direct sum given by $$\h=\mathfrak{t}(R')+\displaystyle\sum_{\alpha\in R'}\lie_{\alpha},$$ where the $\lie_{\alpha}$ are root spaces with $\alpha\in R'$.

\end{itemize}
\end{proposition}

\begin{corollary} Let $\pi:G/K\longrightarrow G/H$ be a homogeneous fibration as above, with totally geodesic fibers isometric to $H/K$ ($\dim H/K\geq 2$) and let $\lie,\sub$ and $\h$ be the Lie algebras of $G,K$ and $H$. Accordingly the current notations, if $\p$ and $\q$ are the vertical and horizontal distribution, i.e, $\p$ is the tangent space to the fiber $H/K$ and $\q$ is the tangent space to the basis space $G/H$, then $\p$ is equal to the $(-B)$-orthogonal direct sum of the root spaces $\lie_{\alpha}$ which $\alpha\in R'$ and hence $\q$ is equal to the $(-B)$-orthogonal direct sum of the root spaces $\lie_{\beta}$ which $\beta\in R\setminus R'$.

\end{corollary}

We can now introduce the construction of the canonical variations $\cvg,\cvh,\cvk,\cvm,\cvn$, of the normal metrics on the maximal flag manifolds $SU(n+1)/T^{n},SO(2n+1)/T^{n},Sp(n)/T^{n}, SO(2n)/T^{n}, G_2/T $, respectively.

% Let $\m=\m_{\alpha_1}\oplus\m_{\alpha_2}\oplus\ldots\oplus\m_r$ be the decomposition of the isotropy representation $\m$ into $\Ad(K)$-invariants inequivalent irreducible submodules $\m_{\alpha_i}$, $\alpha_i\in R^+$, $R$ system of roots of the Cartan subalgebra $\sub$ with relation to $\lie $($G/K$ is a maximal flag manifold). It is known that $\dim \m_{\alpha_i}=2$, for every $i=1,\ldots,r$ and that, since $G/K$ is a maximal flag manifold, each $\m_{\alpha_i}$ is the root space of ${\alpha_i}$. By the properties of the Weyl basis of $\lie$, $[\m_{\alpha},\m_{\beta}]=\m_{\alpha+\beta}$, and $[\m_{\alpha},\m_{\beta}]={0}$ if $\alpha+\beta\notin R$, for any $\alpha,\beta\in R^+$. The Lie algebra $\lie$ decomposes into the $-B$-orthogonal direct sum $$\lie=\h\oplus\q,$$ with the spaces $\sub$ and $\q$  satisfying the relations $$[\h,\h]\subset\h, [\h,\q]\subset\q, [\q,\q]\subset\h,$$ since $G/H$ is a symmetric space and hence $(\h,\q)$ is a symmetric pair.

%Moreover, since $H/K$ is a maximal flag manifold and, by the definition of $\p$, $\h=\sub\oplus\p$, Remark 5 allows us to conclude that $\p=\displaystyle\sum_{\alpha\in R'}\m_{\alpha}$. \cqd
\subsubsection{$(SU(n+1)/T^{n},\cvg), n\geq 2$:}

We observe that for $n=1$ the maximal flag $SU(n+1)/T^{n}$ is $SU(2)/S(U(1)\times U(1))\cong \mathbb{CP}^1$, which is an isotropy irreducible Hermitian symmetric space. 
 
We will denote by $G=SU(n+1)$ the compact simple Lie group whose Lie algebra is $\lie=\mathfrak{su}(n+1)$ and by $T^{n}\subset G$ a maximal torus, given by $T^{n}=S(U(1)\times\ldots\times U(1))=S(U(1)^{n+1})$. The maximal flag manifold associated with $\mathfrak{su}(n+1)$ is $SU(n+1)/T^{n}$. The Lie algebra $\mathfrak{su}(n+1)$ decomposes into the $(-B)$-orthogonal direct sum $$\mathfrak{su}(n+1)=\sub\oplus\m,$$ where $B$ is the Cartan-Killing form of $\mathfrak{su}(n+1)$, defined by $B(X,Y)=2(n+1)\text{Tr}(XY),$ $X,Y\in \mathfrak{su}(n+1),$ $\sub$ is the Lie algebra of $T^{n}$, maximal abelian Lie subalgebra of $\mathfrak{su}(n+1)$ formed by the diagonal matrices of the form $$\sub=\left\{\sqrt{-1}\cdot\text{diag}(a_1,\ldots,a_{n+1});\displaystyle\sum_{i=1}^n{a_i}=0\right\}$$ and $\m$ the $\Ad (T^{n})$-invariant isotropy representation of $T^{n}$. With respect to the Cartan subalgebra $\sub^{\mathbb{C}}$, the root system of $\mathfrak{su}(n+1)^{\mathbb{C}}$ can be chosen as $$R=\{\alpha_{ij}=\pm(\lambda_i-\lambda_j); i\neq j\},$$ where $\lambda_i$ is given by $\text{diag}(a_1,\ldots,a_{n+1})\mapsto \lambda_i(\text{diag}(a_1,\ldots,a_{n+1}))=a_i$, for each $1\leq i\leq n+1$. We fix a system of simple roots to be $\Sigma=\{\alpha_{ii+1};1\leq i\leq n\}$ and with respect to $\Sigma$ the positive roots are given by $$R^+=\{\alpha_{ij}=\lambda_i-\lambda_j; i<j\}.$$In this case one has $\dfrac{n(n+1)}{2}$ positive roots, hence $\m$ decomposes into $\dfrac{n(n+1)}{2}$ pairwise inequivalent irreducible $\Ad(T^{n})$-modules, that is,
\begin{equation}
\m=\displaystyle\bigoplus_{i<j\leq n+1}\m_{\alpha_{ij}}, \label{decompsu}
\end{equation}
where each two dimensional irreducible submodule $\m_{\alpha}$ above defined is generated by $\{A_{\alpha},S_{\alpha}\}$, where $A_{\alpha}=X_{\alpha}+X_{-\alpha}$, $S_{\alpha}=\sqrt{-1}(X_{\alpha}-X_{-\alpha})$ and $X_{\alpha}$ belongs to the Weyl basis of $\mathfrak{su}(n+1).$ Remembering that the root vectors $\{H_{\alpha_{12}},\ldots,H_{\alpha_{nn+1}}\}\cup\{X_{\alpha}\in\mathfrak{su}(n+1)_{\alpha}^{\mathbb{C}}\}$ of this basis satisfy $B(X_{\alpha},X_{-\alpha})=-1$ and $[X_{\alpha},X_{-\alpha}]=-H_{\alpha}\in\sub.$ 

Setting $H=S(U(1)\times U(n))$ and $K=T^{n}$, we have $K\subsetneq H\subsetneq G$, with $G,H$ and $K$ compact connected Lie groups. Consider the canonical map $$\pi:SU(n+1)/T^{n}\longrightarrow SU(n+1)/S(U(1)\times U(n)).$$ Let $\lie =\mathfrak{su}(n+1)$, $\h=\mathfrak{s}(\mathfrak{u}(1)\oplus\mathfrak{u}(n))$ and $\sub$ be the Lie algebras of $G$, $H$ and $K$, respectively. As we saw previously, since $K,H$ and $G$ are compacts, we can consider a $\Ad(H)$-invariant $(-B)$-orthogonal complement $\mathfrak{q}$ to $\h$ in $\lie$, and a $\Ad(K)$-invariant $(-B)$-orthogonal complement $\mathfrak{p}$ to $\sub$ in $\h$.

The Lie algebra $\lie$ decomposes into the sum $$\lie=\mathfrak{su}(n+1)=\sub\oplus \mathfrak{p}\oplus\mathfrak{q}=\sub\oplus\m,$$ hence
\begin{equation}
\pi:(SU(n+1)/T^{n},g)\longrightarrow (SU(n+1)/S(U(1)\times U(n)),\breve{g}) \label{homfibsu}
\end{equation}
is a Riemannian submersion with totally geodesic fibers isometric to $(SU(n)/T^{n-1},\hat{g})$, since $$S(U(1)\times U(n))/T^{n}=S(U(1)\times U(n))/S(U(1)^{n+1})=SU(n)/T^{n-1},$$$T^{n-1}=S(U(1)^{n})$, with $g$ the normal metric determined by the inner product $(-B)|_{\m}$, $\hat{g}$ the metric given by $(-B)|_{\mathfrak{p}}$ and $\breve{g}$ defined by the inner product $(-B)|_{\mathfrak{q}}$.
Since $\m=\mathfrak{p}\oplus\mathfrak{q}$ and $\dim SU(n)/T^{n-1}=n(n-1)$, and due fact that $\dim\m_{\alpha}=2$, for each positive root $\alpha$, we have  the vertical space $\mathfrak{p}$ (tangent to the fiber $SU(n)/T^{n-1}$) is equal to the direct sum of $\dfrac{n(n-1)}{2}$ irreducible isotropy summands of $\m$ while $\mathfrak{q}$ is the direct sum of $\dfrac{n(n+1)}{2}-\dfrac{n(n-1)}{2}=n$ isotropy summands, since the dimension of the total space is equal to $n(n+1)$.

In order to find the submodules that compose these distributions, observe that the fiber $SU(n)/T^{n-1}$ of the original fibration is a maximal flag manifold and the space $\mathfrak{p}$ is exactly the isotropy representation of $T^{n-1}$, which in turn decomposes into $\dfrac{n(n-1)}{2}$ pairwise inequivalent irreducible $\Ad(T^{n-1})$-modules equal to the root spaces relative to the positive roots in $R'\subset R$,
$$R'=\{\alpha_{ij}=\lambda_i-\lambda_j; i<j\leq n\}.$$The Riemannian metric $\hat{g}$ on the fiber $SU(n)/T^{n-1}$, represented by $(-B)|_{\mathfrak{p}}$, is the normal metric represented by the inner product $$g_{eT^{n-1}}=\displaystyle\sum_{i<j\leq n}(-B)|_{\m_{\alpha_{ij}}}$$and the homogeneous metric $\breve{g}$ defined by the inner product $(-B)|_{\mathfrak{q}}$ makes the basis $(SU(n+1)/S(U(1)\times U(n)),\breve{g})$ an isotropy irreducible compact symmetric space.

Hence, by scaling the normal metric $$g_{eT^{n}}=\displaystyle\sum_{i<j\leq n+1}(-B)|_{\m_{\alpha_{ij}}}$$ on the total space $SU(n+1)/T^{n}$ in the direction of the fibers by $t^2$, i.e., multiplying by $t^2, t>0$, the parcels $(-B)|_{\m_{\alpha_{ij}}},1\leq i<j\leq n$ in the expression of $g$ we obtain the canonical variation $\cvg$ of the metric $g$,  according \eqref{gt}.

\subsubsection{$(SO(2n+1)/T^n,\cvh), n\geq 2, n\neq 3$:} 
 
We will denote by $G=SO(2n+1)$ the compact simple Lie group whose Lie algebra is $\lie=\mathfrak{so}(2n+1)$ and by $T^n\subset G$ a maximal torus, given by $T^n=U(1)\times\ldots\times U(1)=U(1)^n$. The maximal flag manifold associated with $\mathfrak{so}(2n+1)$ is $SO(2n+1)/T^n$. The Lie algebra $\mathfrak{so}(2n+1)$ decomposes into the $(-B)$-orthogonal direct sum $$\mathfrak{so}(2n+1)=\sub\oplus\m,$$ where $B$ is the Cartan-Killing form of $\mathfrak{so}(2n+1)$, defined by $B(X,Y)=(2n-1)\text{Tr}(XY),$ $X,Y\in \mathfrak{so}(2n+1),$ $\sub$ is the Lie algebra of $T^n$, maximal abelian Lie subalgebra of $\mathfrak{so}(2n+1)$ formed by the diagonal matrices of the form $$\sub=\left\{\sqrt{-1}\cdot\text{diag}(0,a_1,\ldots,a_n,-a_1,\ldots,-a_n);a_i\in\R\right\}$$ and $\m$ the $\Ad (T^n)$-invariant isotropy representation of $T^n$. With respect the Cartan subalgebra $\sub^{\mathbb{C}}$, the root system of $\mathfrak{so}(2n+1)^{\mathbb{C}}$ can be chosen as by $$R=\{\pm\lambda_i\pm\lambda_j,\pm\lambda_k; 1\leq i<j\leq n,1\leq k\leq n\},$$ where $\lambda_i$ is given by $\text{diag}(a_1,\ldots,a_n)\mapsto \lambda_i(\text{diag}(a_1,\ldots,a_n))=a_i$, for each $1\leq i\leq n$. The system of simple roots is $\Sigma=\{\lambda_1-\lambda_2,\ldots,\lambda_{n-1}-\lambda_n,\lambda_n\}$ and with respect to $\Sigma$ the positive roots are given by $R^+=\{\lambda_i-\lambda_j,\lambda_i+\lambda_j,\lambda_k;1\leq i<j\leq n,1\leq k\leq n\}$. In this case one has $n^2$ positive roots, hence $\m$ decomposes into $n^2$ pairwise inequivalent irreducible $\Ad(T^n)$-modules as follows 
\begin{equation}
\m=\displaystyle\bigoplus_{\alpha\in R^+}\m_{\alpha}, \label{decompso}
\end{equation}
where each two dimensional irreducible submodule $\m_{\alpha}$ above is generated by $\{A_{\alpha},S_{\alpha}\}$, where $A_{\alpha}=X_{\alpha}+X_{-\alpha}$, $S_{\alpha}=\sqrt{-1}(X_{\alpha}-X_{-\alpha})$ and $X_{\alpha}$ belongs to the Weyl basis of $\mathfrak{so}(2n+1)^{\mathbb{C}}.$ 

Setting $H=SO(2n)$ and $K=T^n$, we have that $K\subsetneq H\subsetneq G$, with $G,H$ and $K$ compact connected Lie groups. Consider the canonical map $$\pi:SO(2n+1)/T^n\longrightarrow SO(2n+1)/SO(2n).$$ Let $\lie =\mathfrak{so}(2n+1)$, $\h=\mathfrak{so}(2n)$ and $\sub$ be the Lie algebras of $G$, $H$ and $K$, respectively. As we saw previously, since $K,H$ and $G$ are compacts, we can consider a $\Ad(H)$-invariant $(-B)$-orthogonal complement $\mathfrak{q}$ to $\h$ in $\lie$, and a $\Ad(K)$-invariant $(-B)$-orthogonal complement $\mathfrak{p}$ to $\sub$ in $\h$.

The Lie algebra $\lie$ decomposes into the sum $$\lie=\mathfrak{so}(2n+1)=\sub\oplus \mathfrak{p}\oplus\mathfrak{q}=\sub\oplus\m,$$ hence
\begin{equation}
\pi:(SO(2n+1)/T^n,g)\longrightarrow (SO(2n+1)/SO(2n),\breve{g})\label{homfibso1}
\end{equation}
is a Riemannian submersion with totally geodesic fibers isometric to the maximal flag manifold $(H/K,\hat{g})$, $$H/K=SO(2n)/T^n,$$ where $g$ is the normal metric determined by the inner product $-B|_{\m}$, $\hat{g}$ the metric given by $(-B)|_{\mathfrak{p}}$ and $\breve{g}$ defined by the inner product $(-B)|_{\mathfrak{q}}$.
Since $\m=\mathfrak{p}\oplus\mathfrak{q}$ and $\dim SO(2n)/T^n=2n(n-1)$, and due the fact that $\dim\m_{\alpha}=2$, for each positive root $\alpha$, we have that the vertical space $\mathfrak{p}$ (tangent to the fiber $SO(2n)/T^n$) is equal to the direct sum of $n(n-1)$ irreducible isotropy summands of $\m$ while $\mathfrak{q}$ is the direct sum of $n^2-n(n-1)=n$ isotropy summands, since the  dimension of the total space is $2n^2$.

Since the fiber $SO(2n)/T^n$ of the original fibration is a maximal flag manifold, the space $\mathfrak{p}$ is exactly the isotropy representation of $T^n$, which in turn decomposes into $n(n-1)$ pairwise inequivalent irreducible $\Ad(T^{n})$-modules equal to the root spaces relative to the positive roots in $R'\subset R$,
$$R'=\{\pm\lambda_i\pm\lambda_j; 1\leq i<j\leq n\}.$$The Riemannian metric $\hat{g}$ on the fiber $SO(2n)/T^n$, represented by $(-B)|_{\mathfrak{p}}$, is the normal metric determined by the inner product $$g_{eT^n}=\displaystyle\sum_{\alpha\in R'\cap R^+}(-B)|_{\m_{\alpha}}$$and the homogeneous metric $\breve{g}$ defined by the inner product $(-B)|_{\mathfrak{q}}$ makes the basis $(SU(2n+1)/SO(2n),\breve{g})$ a irreducible Hermitian symmetric space, isometric to the round sphere.

We remark that for $n=2$, we have the fibration $$\pi:(SO(5)/T^2,g)\longrightarrow (SO(5)/SO(4)=S^4,\breve{g})$$with totally geodesic fibers isometric to $(H/K,\hat{g})$, $H/K=SO(4)/T^2\cong SO(4)/SO(2)\times SO(2)\cong S^2\times S^2$, $g$ the normal metric determined by the inner product $(-B)|_{\m}$, $\hat{g}$ the product metric induced by $(-B)|_{\mathfrak{p}}$ and $\breve{g}$ defined by the inner product $(-B)|_{\mathfrak{q}}$.

Hence, by scaling the normal metric $$g_{eT^{n}}=\displaystyle\sum_{\alpha\in R^+}(-B)|_{\m_{\alpha}}$$ on the total space $SO(2n+1)/T^n$ in the direction of the fibers by $t^2$, i.e., multiplying by $t^2, t>0$, all the $(-B)|_{\m_{\alpha}},\alpha\in R'\cap R^+$ that appear in the expression of $g$, we obtain the canonical variation $\cvh$ of $g$.

\subsubsection{$(Sp(n)/T^n,\cvk), n\geq 3$:} 

Let $G=Sp(n)$ be the compact simple Lie group whose Lie algebra is $\lie=\mathfrak{sp}(n)$ and by $T^n\subset G$ a maximal torus there, given by $T^n=U(1)\times\ldots\times U(1)=U(1)^n$. The maximal flag manifold associated is $Sp(n)/T^n$ and the Lie algebra $\mathfrak{sp}(n)$ decomposes into the $(-B)$-orthogonal direct sum $$\mathfrak{sp}(n)=\sub\oplus\m,$$ where $B$ is its the Cartan-Killing form defined by $B(X,Y)=2(n+1)\text{Tr}(XY),$ $X,Y\in \mathfrak{sp}(n),$ $\sub$ is the Lie algebra of $T^n$, maximal abelian Lie subalgebra of $\mathfrak{sp}(n)$ formed by the diagonal matrices of the form $$\sub=\left\{\sqrt{-1}\cdot\text{diag}(a_1,\ldots,a_n,-a_1,\ldots,-a_n);a_i\in\R\right\}$$ and $\m$ the $\Ad (T^n)$-invariant isotropy representation of $T^n$. With respect the Cartan subalgebra $\sub^{\mathbb{C}}$, the root system of $\mathfrak{sp}(n)^{\mathbb{C}}$ can be chosen as $$R=\{\pm\lambda_i\pm\lambda_j,\pm2\lambda_k; 1\leq i<j\leq n,1\leq k\leq n\},$$ where $\lambda_i$ is given by $\text{diag}(a_1,\ldots,a_n)\mapsto \lambda_i(\text{diag}(a_1,\ldots,a_n))=a_i$, for each $1\leq i\leq n$. The system of simple roots is $\Sigma=\{\lambda_1-\lambda_2,\ldots,\lambda_{n-1}-\lambda_n,2\lambda_n\}$ and with respect to $\Sigma$ the positive roots are given by $R^+=\{\lambda_i-\lambda_j,\lambda_i+\lambda_j,2\lambda_k;1\leq i<j\leq n,1\leq k\leq n\}$. In this case one has $n^2$ positive roots, hence $\m$ decomposes into $n^2$ pairwise inequivalent irreducible $\Ad(T^n)$-modules, namely
\begin{equation}
\m=\displaystyle\bigoplus_{\alpha\in R^+}\m_{\alpha}, \label{decompsp}
\end{equation}
where each two dimensional irreducible submodule $\m_{\alpha}$ above is generated by $\{A_{\alpha},S_{\alpha}\}$, where $A_{\alpha}=X_{\alpha}+X_{-\alpha}$, $S_{\alpha}=\sqrt{-1}(X_{\alpha}-X_{-\alpha})$ and $X_{\alpha}$ belongs to the Weyl basis of $\mathfrak{sp}(n)^{\mathbb{C}}.$ 

Setting $H=U(n)$ and $K=T^n$, we have $K\subsetneq H\subsetneq G$, with $G,H$ and $K$ compact connected Lie groups. Consider the canonical map $$\pi:Sp(n)/T^n\longrightarrow Sp(n)/U(n).$$ Let $\h=\mathfrak{u}(n)$ and $\sub$ be the Lie algebra of $H$. Since $K,H$ and $G$ are compacts, we can consider an $\Ad(H)$-invariant $(-B)$-orthogonal complement $\mathfrak{q}$ to $\h$ in $\lie$, and a $\Ad(K)$-invariant $-B$-orthogonal complement $\mathfrak{p}$ to $\sub$ in $\h$.

The Lie algebra $\lie$ decomposes into the sum $$\lie=\mathfrak{sp}(n)=\sub\oplus \mathfrak{p}\oplus\mathfrak{q}=\sub\oplus\m,$$ and follows that 
\begin{equation}
\pi:(Sp(n)/T^n,g)\longrightarrow (Sp(n)/U(n),\breve{g})\label{homfibsp}
\end{equation}
is a Riemannian submersion with totally geodesic fibers isometric to the maximal flag manifold $(H/K,\hat{g})$, $$H/K=U(n)/T^n\cong SU(n)/S(U(1)\times \ldots\times U(1))=SU(n)/S(U(1)^n),$$ $g$ being the normal metric determined by the inner product $(-B)|_{\m}$, $\hat{g}$ the metric given by $(-B)|_{\mathfrak{p}}$ and $\breve{g}$ defined by the inner product $(-B)|_{\mathfrak{q}}$.
Since $\m=\mathfrak{p}\oplus\mathfrak{q}$ and $\dim H/K= \dim U(n)/T^n =n(n-1)$,  and due the fact that $\dim\m_{\alpha}=2$ for each positive root $\alpha$, we have that the vertical space $\mathfrak{p}$ (tangent to the fiber $SU(n)/S(U(1)^n)$) is equal to the direct sum of $\dfrac{n(n-1)}{2}$ irreducible isotropy summands of $\m$ while $\mathfrak{q}$ is the direct sum of $n^2-\dfrac{n(n-1)}{2}=\dfrac{n(n+1)}{2}$ isotropy summands, since the dimension of the total space is $2n^2$.

Since the fiber $SU(n)/S(U(1)^n)$ of the original fibration is a maximal flag manifold, the space $\mathfrak{p}$ is exactly the isotropy representation of the maximal torus $T^{n-1}_{SU(n)}=S(U(1)^n)$ of $SU(n)$, which in turn decomposes into $n(n-1)$ pairwise inequivalent irreducible $\Ad(T^{n-1}_{SU(n)})$-modules equal to the root spaces relative to the positive roots in $R'\subset R$,
$$R'=\{\pm(\lambda_i-\lambda_j); 1\leq i<j\leq n\}.$$The Riemannian metric $\hat{g}$ on the fiber $SU(n)/S(U(1)^n)$, represented by $(-B)|_{\mathfrak{p}}$, is the normal metric rerpesented by the inner product $$g_{eT^{n-1}_{SU(n)}}=\displaystyle\sum_{\alpha\in R'\cap R^+}(-B)|_{\m_{\alpha}}$$and the homogeneous metric $\breve{g}$ defined by the inner product $(-B)|_{\mathfrak{q}}$ makes the basis $(Sp(n)/U(n),\breve{g})$ a irreducible Hermitian symmetric space.

By scaling the normal metric $$g_{eT^{n}}=\displaystyle\sum_{\alpha\in R^+}(-B)|_{\m_{\alpha}}$$ on the total space $Sp(n)/T^n$ in the direction of the fibers by $t^2$, i.e., multiplying by $t^2, t>0$, all the $(-B)|_{\m_{\alpha}},\alpha\in R'\cap R^+$ in the expression of $g$, we obtain the canonical variation $\cvk$ of $g$.

\subsubsection{$(SO(2n)/T^n,\cvm), n\geq 4$:} 

Let $G=SO(2n)$ be the compact simple Lie group whose Lie algebra is $\lie=\mathfrak{so}(2n)$ and by $T^n\subset G$ a maximal torus with $T^n=U(1)\times\ldots\times U(1)=U(1)^n$. The maximal flag manifold associated is $SO(2n)/T^n$ and the Lie algebra $\mathfrak{so}(2n)$ decomposes into the $(-B)$-orthogonal direct sum $$\mathfrak{so}(2n)=\sub\oplus\m,$$ where $B$ is its the Cartan-Killing form defined by $B(X,Y)=2(n-1)\text{Tr}(XY),$ $X,Y\in \mathfrak{sp}(n),$ $\sub$ is the Lie algebra of $T^n$, maximal abelian Lie subalgebra of $\mathfrak{so}(2n)$ formed by the diagonal matrices of the form $$\sub=\left\{\sqrt{-1}\cdot\text{diag}(a_1,\ldots,a_n,-a_1,\ldots,-a_n);a_i\in\R\right\}$$ and $\m$ the $\Ad (T^n)$-invariant isotropy representation of $T^n$. With respect the Cartan  subalgebra $\sub^{\mathbb{C}}$, the root system of $\mathfrak{so}(2n)^{\mathbb{C}}$ can be chosen as $$R=\{\pm\lambda_i\pm\lambda_j; 1\leq i<j\leq n,\},$$ according to the above notation. The system of simple roots is $$\Sigma=\{\lambda_1-\lambda_2,\ldots,\lambda_{n-1}-\lambda_n,\lambda_{n-1}+\lambda_n\}$$and with respect to $\Sigma$ the positive roots belong to $$R^+=\{\lambda_i-\lambda_j,\lambda_i+\lambda_j;1\leq i<j\leq n\}.$$ One has $n(n-1)$ positive roots, hence $\m$ decomposes into $n(n-1)$ pairwise inequivalent irreducible $\Ad(T^n)$-modules, namely 
\begin{equation}
\m=\displaystyle\bigoplus_{\alpha\in R^+}\m_{\alpha}, \label{decompso2}
\end{equation}
where each two dimensional irreducible submodule $\m_{\alpha}$ above is generated by $\{A_{\alpha},S_{\alpha}\}$, where $A_{\alpha}=X_{\alpha}+X_{-\alpha}$, $S_{\alpha}=\sqrt{-1}(X_{\alpha}-X_{-\alpha})$ and $X_{\alpha}$ belongs to the Weyl basis of $\mathfrak{so}(2n)^{\mathbb{C}}.$ 

Take the canonical map $$\pi:SO(2n)/T^n\longrightarrow SO(2n)/U(n).$$Let $\h=\mathfrak{u}(n)$ and $\sub$ be the Lie algebras $H$ and $K=T^n$, respectively. Since $K,H$ and $G$ are compacts, we can consider a $\Ad(H)$-invariant $(-B)$-orthogonal complement $\mathfrak{q}$ to $\h$ in $\lie$, and a $\Ad(K)$-invariant $(-B)$-orthogonal complement $\mathfrak{p}$ to $\sub$ in $\h$.

The Lie algebra $\lie$ decomposes into the sum $$\lie=\mathfrak{so}(2n)=\sub\oplus \mathfrak{p}\oplus\mathfrak{q}=\sub\oplus\m,$$ and follows that 
\begin{equation}
\pi:(SO(2n)/T^n,g)\longrightarrow (SO(n)/U(n),\breve{g})\label{homfibso2}
\end{equation}
is a Riemannian submersion with totally geodesic fibers isometric to the maximal flag manifold $(H/K,\hat{g})$, $$H/K=U(n)/T^n\cong SU(n)/S(U(1)\times \ldots\times U(1))=SU(n)/S(U(1)^n),$$ $g$ being the normal metric determined by the inner product $(-B)|_{\m}$, $\hat{g}$ the metric given by $(-B)|_{\mathfrak{p}}$ and $\breve{g}$ defined by the inner product $(-B)|_{\mathfrak{q}}$.
Since $\m=\mathfrak{p}\oplus\mathfrak{q}$ and $\dim H/K=\dim U(n)/T^n =n(n-1)$, and due the fact that $\dim\m_{\alpha}=2$, for each positive root $\alpha$, we have that the vertical space $\mathfrak{p}$ (tangent to the fiber $SU(n)/S(U(1)^n)$) is equal to the direct sum of $\dfrac{n(n-1)}{2}$ irreducible isotropy summands of $\m$ while $\mathfrak{q}$, from the fact that $2n(n-1)$ is the dimension of the total space, is equal to the direct sum of $n(n-1)-\dfrac{n(n-1)}{2}=\dfrac{n(n-1)}{2}$ isotropy summands. 

Since the fiber $SU(n)/S(U(1)^n)$ of the original fibration is a maximal flag manifold, the space $\mathfrak{p}$ is exactly the isotropy representation of $T^{n-1}_{SU(n)}=S(U(1)^n)$, which in turn decomposes into $n(n-1)$ pairwise inequivalent irreducible $\Ad(T^{n-1}_{SU(n)})$-modules equal to the root spaces relative to the positive roots in $R'\subset R$,
$$R'=\{\pm(\lambda_i-\lambda_j); 1\leq i<j\leq n\}.$$The Riemannian metric $\hat{g}$ on the fiber $SU(n)/S(U(1)^n)$, represented by $(-B)|_{\mathfrak{p}}$, is the normal metric given by the inner product $$g_{eT^{n-1}_{SU(n)}}=\displaystyle\sum_{\alpha\in R' \cap R^+}(-B)|_{\m_{\alpha}}$$and the homogeneous metric $\breve{g}$ on the basis, defined by the inner product $(-B|)_{\mathfrak{q}}$, makes $(Sp(n)/U(n),\breve{g})$ a irreducible Hermitian symmetric space.

By scaling the normal metric $$g_{eT^{n}}=\displaystyle\sum_{\alpha\in R^+}(-B)|_{\m_{\alpha}}$$ on the total space $SO(2n)/T^n$ in the direction of the fibers by $t^2$, i.e., multiplying by $t^2, t>0$, all the $(-B)|_{\m_{\alpha}},\alpha\in R'\cap R^+$ in the expression of $g$, we obtain the canonical variation $\cvm$ of $g$.

\subsubsection{$(G_2/T,\cvn)$}\label{G2cv}

Let $\lie_2$ be the exceptional complex simple Lie algebra of type $G_2$. We will denote by $G_2$ the compact simple Lie group whose Lie algebra is $\lie_2$ and by $T\subset G_2$ a maximal torus there, where $T=U(1)\times U(1)$. The full flag manifold associated with $\lie_2$ is $G_2/T$. The Lie algebra $\lie_2$ decomposes into the $(-B)$-orthogonal direct sum $$\lie_2=\sub\oplus\m,$$ where $B$ is the Cartan-Killing form of $G_2$, defined by $B(X,Y)=\text{Tr}(\text{ad}(X)\text{ad}(Y)),$ $X,Y\in\lie_2,$ $\sub$ is the Lie algebra of $T$, maximal abelian Lie subalgebra of $\lie_2$ formed by the diagonal traceless matrices belonging to $\lie_2$ and $\m$ its isotropy representation of $T$. With respect the Cartan subalgebra $\sub$, the root system of $\lie_2$ can be chosen as by $$\{\pm\alpha_1,\pm\alpha_2,\pm(\alpha_1+\alpha_2),\pm(\alpha_1+2\alpha_2),\pm(\alpha_1+3\alpha_2),\pm(2\alpha_1+3\alpha_2)\},$$and we fix a system of simple roots to be $\Sigma=\{\alpha_1,\alpha_2\}.$ With respect to $\Pi$ the positive roots are given by $$\{\alpha_1,\alpha_2,(\alpha_1+\alpha_2),(\alpha_1+2\alpha_2),(\alpha_1+3\alpha_2),(2\alpha_1+3\alpha_2)\}.$$The maximal root is $\alpha=2\alpha_1+3\alpha_2.$ The angle between $\alpha_1$ and $\alpha_2$ is $5\pi/6$ and we have $\left\|\alpha_1\right\|=\sqrt{3}\left\|\alpha_2\right\|$. Moreover, the roots of $\lie_2$ form successive angles of $\pi/6$. Since in this case one has six positive roots, the space $\m$ decomposes into six pairwise inequivalent irreducible $\Ad(T)$-modules as follows $$\m=\m_{\alpha_1}\oplus\m_{\alpha_2}\oplus\m_{\alpha_1+\alpha_2}\oplus\m_{\alpha_1+2\alpha_2}\oplus\m_{\alpha_1+3\alpha_2}\oplus\m_{2\alpha_1+3\alpha_2},$$ where each two dimensional irreducible submodule $\m_{\alpha}$ above is generated by $\{A_{\alpha},S_{\alpha}\}$, where $A_{\alpha}=X_{\alpha}+X_{-\alpha}$, $S_{\alpha}=\sqrt{-1}(X_{\alpha}-X_{-\alpha})$ and $X_{\alpha}$ belongs to the Weyl basis of $\lie_2$.

Setting $G=G_2, H=SO(4)\cong SO(3)\times SO(3)$ and $K=T=U(1)\times U(1)\cong SO(2)\times SO(2)$, we have that $K\subsetneq H\subsetneq G$, with $G,H$ and $K$ compact connected Lie groups. Consider the canonical map $$\pi:G_2/T\longrightarrow G_2/SO(4).$$ Let $\lie =\lie_2$, $\h=\mathfrak{so}(4)\cong \mathfrak{su}(2)\oplus \mathfrak{su}(2)$ and $$\sub=\left\{d=\text{diag}(ia,ib,-i(a+b))\in\mathfrak{sl}(3)\subset\lie_2; a,b\in\R\right\}$$be the Lie algebras of $G$, $H$ and $K$, respectively. As we saw previously, since $K,H$ and $G$ are compacts, we can consider a $\Ad(H)$-invariant $(-B)$-orthogonal complement $\mathfrak{q}$ to $\h=\mathfrak{so}(4)\cong \mathfrak{su}(2)\oplus \mathfrak{su}(2)$ in $\lie_2$, and a $\Ad(K)$-invariant $(-B)$-orthogonal complement $\mathfrak{p}$ to $\sub$ in $\h$.

The Lie algebra $\lie_2$ decomposes into the sum $$\lie_2=\sub\oplus \mathfrak{p}\oplus\mathfrak{q}=\sub\oplus\m,$$ therefore $$\pi:(G_2/T,g)\longrightarrow (G_2/SO(4),\breve{g})$$ is a Riemannian submersion with totally geodesic fibers isometric to $(H/K,\hat{g})$, $H/K=SO(4)/T\cong SO(4)/SO(2)\times SO(2)\cong S^2\times S^2$, $g$ the normal metric determined by the inner product $(-B)|_{\m}$, $\hat{g}$ the metric given by $(-B)|_{\mathfrak{p}}$ and $\breve{g}$ defined by the inner product $(-B)|_{\mathfrak{q}}$.

Since $\m=\mathfrak{p}\oplus\mathfrak{q}$ and $\dim SO(4)/T\cong S^2\times S^2=4$, and due to the fact that $\dim\m_{\alpha}=2$ for each positive root $\alpha$, we have that $\mathfrak{q}$ is equal to the direct sum of four irreducible submodules of $\m$ while $\mathfrak{p}$ is the direct sum of two irreducible submodules. In order to find the submodules that composes each of these distributions, we apply the relations $[\h,\mathfrak{q}]\subset\mathfrak{q}$ and $[\sub,\mathfrak{p}]\subset\mathfrak{p}$ and the same property of the Weyl basis in the previous example, obtaining $\mathfrak{p}=\m_{\alpha_1+\alpha_2}\oplus\m_{\alpha_1+3\alpha_2}$ and $\mathfrak{q}=\m_{\alpha_1}\oplus\m_{\alpha_2}\oplus\m_{\alpha_1+2\alpha_2}\oplus\m_{2\alpha_1+3\alpha_2}$, since $\h=\sub\oplus\mathfrak{p}$. Hence, by scaling the normal metric 
\begin{eqnarray*}
g_{eT}&=&(-B)|_{\m_{\alpha_1}}+(-B)|_{\m_{\alpha_2}}+(-B)|_{\m_{\alpha_1+\alpha_2}}+(-B)|_{\m_{\alpha_1+2\alpha_2}}\\
      & &+(-B)|_{\m_{\alpha_1+3\alpha_2}}+(-B)|_{\m_{2\alpha_1+3\alpha_2}},
\end{eqnarray*}
of the total space $G_2/T$ in the direction of the fibers, i.e., multiplying by $t^2, t>0$, the component $(-B)|_{\m_{\alpha_1+\alpha_2}}+(-B)|_{\m_{\alpha_1+3\alpha_2}}$ in the expression of $g$ we obtain the following canonical variation $\cvn$ of the metric $g$
\begin{eqnarray*}
(\cvn)_{eT}&=&(-B)|_{\m_{\alpha_1}}+(-B)|_{\m_{\alpha_2}}+t^2(-B)|_{\m_{\alpha_1+\alpha_2}}+(-B)|_{\m_{\alpha_1+2\alpha_2}}\\
           & & +t^2(-B)|_{\m_{\alpha_1+3\alpha_2}}\linebreak +(-B)|_{\m_{2\alpha_1+3\alpha_2}},
\end{eqnarray*}
according \eqref{gt}.

\subsection{Spectra of Maximal Flag Manifolds and Symmetric spaces} \label{spectraspaces}

In \cite{yamaguchi}, S. Yamaguchi describes the spectrum of the Laplacian defined on $C^{\infty}(G/T,g)$, where $G/T$ is a maximal flag manifold equipped with a normal homogeneous metric $g$. In this case, we have the following description of $\sigma(\Delta_g)$. 

Let $\lie$ be the Lie algebra of $G$, $\lie$ complex and simple Lie algebra. Denote by $\mathfrak{h}$ the Cartan subalgebra of $\lie$, $R$ the root system of $(\lie,\mathfrak{h})$ and by $\Sigma=\{\alpha_1,\ldots,\alpha_l\}$ the associated system of simple roots. Let $(\cdot,\cdot)=-B$ the inner product on $\lie$ induced by the Cartan-Killing form $B$. Consider the {\it system of fundamental weights} $\{\omega_1,\cdots,\omega_l\}$ of $\h$ and denote by $\mathcal{P}$ the set of all {integral dominant weights}, $$\mathcal{P}=\{\Lambda=\displaystyle\sum_{i=1}^{l}s_i\omega_i\in\mathfrak{h}^{\ast};s_i\geq 0, s_i\in\mathbb{Z}\}.$$ Let $\widehat{\mathcal{P}}$ the set of all elements of $\mathcal{P}$ which are of class one relative to $\mathfrak{h}$, i.e., the irreducible representation $(\xi_{\Lambda},V_{\Lambda})$ of $\lie$ in $V_{\Lambda}$, with the highest weight equal to $\Lambda\in\mathcal{P}$, has a non zero $\xi_{\Lambda}(\mathfrak{h})$-invariant vector in the representation space $V_{\Lambda}$.  

\begin{theorem}[H. Freudenthal \cite{yamaguchi}] Let $(\xi_{\Lambda},V_{\Lambda})$ be the irreducible representation with highest weight $\Lambda\in\mathcal{P}$ ($\Lambda\neq 0$). Then, $\Lambda\in\widehat{\mathcal{P}}$ if and only if $\Lambda=\displaystyle\sum_{i=1}^{l}p_i\alpha_i, p_i\geq 1$, $p_i\in\mathbb{Z}$, $1\leq i\leq l$.
\end{theorem}

We identify $\Lambda=\displaystyle\sum_{i=1}^{l}p_i\alpha_i\in\widehat{\mathcal{P}}$ with $P=(p_1,\ldots,p_l)$.

\begin{theorem}[S. Yamaguchi \cite{yamaguchi}] \label{spectra} If $G/T$ is a maximal flag manifold associated with a complex classical simple Lie algebra $\lie$ of the type $A_n,B_n,C_n,D_n$ or with the exceptional Lie algebra of the type $G_2$, the spectrum of the Laplacian $\Delta_g$, $g$ normal metric induced by the Cartan-Killing form $B$ of $\lie$, is determined by the irreducible representations with highest weight $\Lambda\in\widehat{\mathcal{P}}$, in such a way that each eigenvalue $\mu\in\sigma(\Delta_g)$ is given by $\mu(\Lambda)$, for some irreducible representation $(\xi_{\Lambda},V_{\Lambda})$ with highest weight $\Lambda=\displaystyle\sum_{i=1}^{l}p_i\alpha_i$, as follows:
\begin{flushleft}
\bf{\text{Type} $A_n,n\geq 1$:}
\end{flushleft}

\begin{enumerate}
\item[\textbf{(1)}] $\mu(\Lambda)=\dfrac{1}{n+1}\left\{\displaystyle\sum_{i=1}^{n}p_i^2-\displaystyle\sum_{i=1}^{n-1}p_ip_{i+1}+\displaystyle\sum_{i=1}^{n}p_i\right\},
$
\item[\textbf{(2)}]$\left\{
  \begin{array}{rcl} 
  2p_1-p_2&\geq& 0\\
  -p_1+2p_2-p_3&\geq&0\\
  \vdots& & \\
	-p_{n-2}+2p_{n-1}-p_n&\geq&0\\
	-p_{n-1}+2p_n&\geq& 0\\
  \end{array}
  \right.,$ 
  
\item[\textbf{(3)}]\text{First positive eigenvalue} $\mu_1=\mu(1,\ldots,1)=1,$ $$P_0=(1,\ldots,1)\longleftrightarrow \Lambda_0=\alpha_1+\ldots+\alpha_n \ \ \textbf{highest root},$$ $\xi_{\Lambda_0}$=\text{adjoint representation}.
\end{enumerate}
  
\begin{flushleft}
\bf{\text{Type} $B_n,n\geq 2$:}
\end{flushleft} 

\begin{enumerate}
\item[\textbf{(1)}] $\mu(\Lambda)=\dfrac{1}{4n-2}\left\{\displaystyle\sum_{i=1}^{n-1}2p_i^2+p_n^2-2\displaystyle\sum_{i=1}^{n-1}p_ip_{i+1}+2\displaystyle\sum_{i=1}^{n-1}p_i+p_n\right\},
$
\item[\textbf{(2)}]$\left\{
  \begin{array}{rcl} 
  2p_1-p_2&\geq& 0\\
  -p_1+2p_2-p_3&\geq&0\\
  \vdots& & \\
	p_{n-2}+2p_{n-1}-p_n&\geq&0\\
	-p_{n-1}+p_n&\geq& 0\\
  \end{array}
  \right.,$ 
  
\item[\textbf{(3)}]\text{First positive eigenvalue} $\mu_1=\mu(1,\ldots,1)=\dfrac{n}{2n-1},$ $$P_0=(1,\ldots,1)\longleftrightarrow \Lambda_0=\omega_1=\alpha_1+\ldots+\alpha_n,$$ $\xi_{\Lambda_0}$=\text{representation with highest weight} $\omega_1$.
\end{enumerate}

\begin{flushleft}
\bf{\text{Type} $C_n,n\geq 3$:}
\end{flushleft}

\begin{enumerate}
\item[\textbf{(1)}] $\mu(\Lambda)=\dfrac{1}{2(n+1)}\left\{\displaystyle\sum_{i=1}^{n-1}p_i^2+2p_n^2-\displaystyle\sum_{i=1}^{n-2}p_ip_{i+1}-p_{n-1}p_n+\displaystyle\sum_{i=1}^{n-1}p_i+2p_n\right\},
$
\item[\textbf{(2)}]$\left\{
  \begin{array}{rcl} 
  2p_1-p_2&\geq& 0\\
  -p_1+2p_2-p_3&\geq&0\\
  \vdots& & \\
	-p_{n-3}+2p_{n-2}-p_{n-1}&\geq&0\\
	-p_{n-2}+2p_{n-1}-2p_n&\geq& 0\\
	-p_{n-1}+2p_n&\geq& 0\\
  \end{array}
  \right.,$ 
  
\item[\textbf{(3)}]\text{First positive eigenvalue} $\mu_1=\mu(1,2,\ldots,2,1)=\dfrac{4n-1}{4(n+1)},$ $$P_0=(1,\ldots,1)\longleftrightarrow \Lambda_0=\omega_2,$$ $\xi_{\Lambda_0}$=\text{representation with highest weight} $\omega_2$.
\end{enumerate}

\begin{flushleft}
\bf{\text{Type} $D_n,n\geq 4$:}
\end{flushleft}

\begin{enumerate}
\item[\textbf{(1)}] $\mu(\Lambda)=\dfrac{1}{2n-1}\left\{\displaystyle\sum_{i=1}^np_i^2-\displaystyle\sum_{i=1}^{n-2}p_ip_{i+1}-p_{n-2}p_n+\displaystyle\sum_{i=1}^{n}p_i\right\},
$
\item[\textbf{(2)}]$\left\{
  \begin{array}{rcl} 
  2p_1-p_2&\geq& 0\\
  -p_1+2p_2-p_3&\geq&0\\
  \vdots& & \\
	-p_{n-4}+2p_{n-3}-p_{n-2}&\geq&0\\
	-p_{n-3}+2p_{n-2}-p_{n-1}-p_n&\geq& 0\\
	-p_{n-2}+2p_{n-1}&\geq& 0\\
	-p_{n-2}+2p_n&\geq& 0\\
  \end{array}
  \right.,$ 
  
\item[\textbf{(3)}]\text{First positive eigenvalue} $\mu_1=\mu(1,2,\ldots,2,1,1)=1,$ $$P_0=(1,2,\ldots,2,1,1)\longleftrightarrow \Lambda_0 \ \ \textbf{highest root},$$ $\xi_{\Lambda_0}$=\text{adjoint representation}.
\end{enumerate}

\begin{flushleft}
\bf{\text{Type} $G_2$:}
\end{flushleft}

\begin{enumerate}
\item[\textbf{(1)}] $\mu(\Lambda)=\dfrac{1}{12}\left\{p_1^2 + 3p_2^2-3p_1p_2+p_1+3p_2\right\},
$
\item[\textbf{(2)}]$\left\{
  \begin{array}{rcl} 
  2p_1-3p_2&\geq& 0\\
  -p_1+2p_2&\geq&0\\
  \end{array}
  \right.,$ 
  
\item[\textbf{(3)}]\text{First positive eigenvalue} $\mu_1=\mu(2,1)=\dfrac{1}{2}, P_0=(2,1)\longleftrightarrow \Lambda_0=\omega_1$, \\ $\xi_{\Lambda_0}$=\text{representation with highest weight $\omega_1$}.
\end{enumerate}

\end{theorem}

We note that the base spaces $(G/H,\breve{g})$ of the homogeneous fibrations that have been studied so far in this work are irreducible Hermitian symmetric space of compact type. The spectrum of such spaces are well known, and can be determined as follows (see \cite{Urakawa}, page 64).

Let $G$ be a compact simply connected simple Lie group, $H$
closed subgroup of $G$. Let $\lie,\h$ the Lie algebras of $G$ and $H$, respectively, and $\lie=\h\oplus\q$, the Cartan decomposition. The inner product on $\q$ is $(-B)|_{\q}$, where $B$ is the Cartan-Killing form of $\lie$. Let $\breve{g}$ the $G$-invariant Riemannian metric on $G/H$ induced by $B$. Then, it is known that the spectrum of the Laplacian of $(G/H,\breve{g})$ is given by 
\begin{equation}
\mu(\Lambda)=-B(\Lambda+2\delta,\Lambda),\label{specsym}
\end{equation} with multiplicities 
\begin{equation}
d_{\Lambda}=\prod_{\alpha\in R^{+}}\dfrac{-B(\Lambda+\delta,\alpha)}{-B(\delta,\alpha)}.\label{multsym}
\end{equation}
Here, $\Lambda$ varies over the set $D(G,H)$ of the highest weights of all spherical representations of $(G,H)$, $\delta$ is equal to the sum of the positive roots $\alpha\in R^{+}$ of the complexification $\lie^{\mathbb{C}}$ of $\lie$ relative to the maximal abelian subalgebra $\sub^{\mathbb{C}}$ of $\lie^{\mathbb{C}}$ and $d_{\Lambda}$ is the dimension of the irreducible spherical representation of $(G,H)$ with highest weight $\Lambda.$

For simple compact connected Lie group $G$, Krämer in \cite{Kramer} provides a classification of possible subgroups $H$ along which a set of dominant weights whose integral non-negative combinations give all spherical representations.

In particular, in Krämer's classification the basis of the representations in $D(G_2,SO(4))$ is given by $$\mathcal{B}=\{2\pi_1,2\pi_2\},$$ with $\{\pi_1,\pi_2\}$ being the set of fundamental weights of the maximal abelian Lie algebra $\text{Lie}(T)=\sub$. Therefore, by applying \eqref{specsym} and the fact that all $\Lambda\in D(G_2,SO(4))$ is linear combination of the weights in $\mathcal{B}$, follows the case of the spectrum of the Laplacian on the isotropy irreducible symmetric space $G_2/SO(4)$ equipped with the metric determined by the Cartan-Killing form of $G_2$.

\begin{proposition}\label{spectrumG2} The spectrum of the Laplacian of the isotropy irreducible Hermitian space $(G_2/SO(4),\breve{g})$ is given by $$\sigma(\Delta_{\breve{g}})=\left\{\dfrac{1}{6}(9r+6r^2+5s+6rs+2s^2);\mathbb{Z}\ni r,s\geq 0\right\},$$ with first positive eigenvalue $\beta_1=\dfrac{7}{6}.$
\end{proposition} 

Assume that $B$ is the Cartan-Killing form of $G$, $G/T$ total space of the homogeneous fibration $$\pi:(G/T,g)\longrightarrow (G/H,\breve{g})$$and $\breve{g}$ is the homogeneous metric given by the inner product $(-B)|_{\q}$, $\q$ being a $(-B)$-orthogonal complement to $\text{Lie}(H)=\h$ in $\text{Lie}(G)=\lie$, i.e, $\q$ is the horizontal distribution relative to the submersion $\pi$. Moreover, each $(G/H,\breve{g})$ below is an isotropy irreducible Hermitian symmetric space. 

We will now describe the spectrum of the Laplacian on each basis space of the homogeneous fibrations constructed above on the maximal flag manifolds $SU(n+1)/T^n, SO(2n+1)/T^n, Sp(n)/T^n, SO(2n)/T^n$ and $G_2/T$, respectively.
 
\subsubsection{Spectrum of the Complex Projective Space, $\sigma(\Delta_{\mathbb{CP}^n})$}
The basis space $(G/H,\breve{g})=(SU(n+1)/S(U(1)\times U(n)),\breve{g})$ of the homogeneous fibration $$\pi:(SU(n+1)/T^{n},g)\longrightarrow (SU(n+1)/S(U(1)\times U(n)),\breve{g}),n\geq 2,$$ is the homogeneous realization of the complex projective space $\mathbb{CP}^n$ equipped with a metric homothetical to the Fubini-Study one. In fact, the induced metric $(-B)|_{\q}$ on $\mathbb{CP}^n$, $B$ being the Cartan-Killing form given by $B(X,Y)=2(n+1)\text{Tr}(XY),X,Y\in\mathfrak{su}(n+1)$, is simply $n+1$ times the Fubini-Study metric. According to \cite{taniguchi}, the spectra of the Laplacian acting on functions on the complex projective space $\mathbb{CP}^n=SU(n+1)/S(U(1)\times U(n))$ equipped with the Fubini-Study metric is $$\sigma(\Delta_{FS})=\{\xi_k=k(k+n);k\in\mathbb{N}\}.$$ Hence, the spectrum of the Laplacian on $\mathbb{CP}^n=SU(n+1)/S(U(1)\times U(n))$ with the metric $\breve{g}$ represented by the inner product $-B|_{\q}$ is
\begin{equation}
 \sigma(\Delta_{\breve{g}})=\left\{\beta_k=\dfrac{\xi_k}{n+1}=\dfrac{k(k+n)}{n+1};k\in\mathbb{N}\right\}.\label{specAn}
\end{equation} Note that the first positive eigenvalue in this case is $\beta_1=1.$

\subsubsection{Spectrum of the Round Sphere, $\sigma(\Delta_{S^{2n}})$}

The basis space $(G/H,\breve{g})=(SO(2n+1)/SO(2n)),\breve{g})$ of the homogeneous fibration $$\pi:(SO(2n+1)/T^{n},g)\longrightarrow (SO(2n+1)/SO(2n),\breve{g}), n\geq 2,$$ is the homogeneous realization of the round sphere $S^{2n}$ equipped with a metric homothetical  the canonical metric. In fact, the metric induced by the Cartan-Killing form given by $$B(X,Y)=(2n-1)\text{Tr}(XY),X,Y\in\mathfrak{so}(2n+1),$$ is $2(2n-1)$ times the usual one. According \cite{taniguchi}, the spectra of the Laplacian acting on functions on sphere $S^{2n}=SO(2n+1)/SO(2n)$ with the usual metric $h$ is $$\sigma(\Delta_h)=\{\xi_k=k(k+n-1);k\in\mathbb{N}\}.$$ Hence, the spectrum of the Laplacian on $S^{2n}=SO(2n+1)/SO(2n)$equipped with the metric $\breve{g}$ represented by the inner product $-B|_{\q}$ is
\begin{equation}
 \sigma(\Delta_{\breve{g}})=\{\beta_k=\dfrac{\xi_k}{2(2n-1)}=\dfrac{k(k+2n-1)}{2(2n-1)};k\in\mathbb{N}\}.\label{specBn}
\end{equation} Note that the first positive eigenvalue in this case is $\beta_1=\dfrac{n}{2n-1}.$

\subsubsection{Spectrum of the Symmetric Space $Sp(n)/U(n)$, $\sigma(\Delta_{Sp(n)/U(n)})$}

In the case of the basis space $(G/H,\breve{g})=(Sp(n)/U(n),\breve{g})$ of the homogeneous fibration $\pi:(Sp(n)/T^{n},g)\longrightarrow (Sp(n)/U(n),\breve{g}), n\geq 3,$ Krämer's classification \cite{Kramer} give us the basis of the representations in the set $D(Sp(n),U(n))$ of the highest weights of all spherical representations of the pair $(Sp(n),U(n))$, namely $$\mathcal{B}=\{2\pi_l;1\leq l\leq n\},$$ with $\{\pi_l;1\leq l\leq n\}$ being the fundamental weights of the maximal abelian Lie algebra $\text{Lie}(T^n)=\sub^{\mathbb{C}}$. Therefore, by applying \eqref{specsym} and the fact that all $\Lambda\in D(Sp(n),U(n))$ is a linear combination of the weights in $\mathcal{B}$, we have the spectrum of the Laplacian on the isotropy irreducible symmetric space $Sp(n)/U(n)$ equipped with the metric determined by the Cartan-Killing form $B$ of $Sp(n)$, given by $$B(X,Y)=2(n+1)\text{Tr}(XY),X,Y\in\mathfrak{sp}(n).$$
The first positive eigenvalue in this case is $\beta_1=1$, see \cite{Urakawa}.

\subsubsection{Spectrum of the Symmetric Space $SO(2n)/U(n)$, $\sigma(\Delta_{SO(2n)/U(n)})$}

In the case of the basis space $(G/H,\breve{g})=(SO(2n)/U(n),\breve{g})$ of the homogeneous fibration $$\pi:(SO(2n)/T^{n},g)\longrightarrow (SO(2n)/U(n),\breve{g}), n\geq 4,$$ Krämer's classification \cite{Kramer} give us the basis of the representations in the set $D(SO(2n),U(n))$ of the highest weights of all spherical representations of the pair $(SO(2n),U(n))$, namely $$\mathcal{B}_1=\{\pi_2,\pi_4,\ldots,\pi_{n-2},2\pi_n\},\ \ \text{if} \ \ n \ \ \text{is even}$$and $$\mathcal{B}_2=\{\pi_2,\pi_4,\ldots,\pi_{n-3},\pi_{n-1}+\pi_n\},\ \ \text{if} \ \ n \ \ \text{is odd}$$with $\{\pi_l;1\leq l\leq n\}$ being the fundamental weights of the maximal abelian Lie algebra $\text{Lie}(T^n)=\sub^{\mathbb{C}}$. Therefore, by applying \eqref{specsym} and the fact that all $\Lambda\in D(SO(2n),U(n))$ is linear combination of the weights in $\mathcal{B}_1$ or $\mathcal{B}_2$, according to $n$ is even or odd, we have the spectrum of the Laplacian on the isotropy irreducible symmetric space $SO(2n)/U(n)$ equipped with the metric determined by the Cartan-Killing form $B$ of $SO(2n)$, given by $$B(X,Y)=2(n-1)\text{Tr}(XY),X,Y\in\mathfrak{so}(2n).$$
The first positive eigenvalue in both of the cases, i.e, when $n$ is even or odd is $\beta_1=1,$ see \cite{Urakawa}.

From the above, we obtain the following useful properties of the first positive eigenvalue $\lambda_1(t)$ of the Laplacian $\Delta_t=\Delta_{g_t}$, where $g_t$ is a canonical variation of the normal metrics on the flag manifolds that appear in Theorem \ref{spectra}.

\begin{proposition}\label{primeiroaut} Considering the canonical variations $(SU(n+1)/T^{n},\cvg),\linebreak (SO(2n+1)/T^{n},\cvh), (Sp(n)/T^{n},\cvk),(SO(2n)/T^{n},\cvm)$ and $(G_2/T,\cvn)$, one has the following estimates for the first positive eigenvalue $\lambda_1(t)$ of the Laplacian $\Delta_t$ on the canonical variations $(G/T,g_t)$:

\begin{table}[htb]
\centering

\begin{tabular}{lr}
\hline
$(G/T,g_t)$ & $\lambda_1(t), 0<t\leq 1$\\
\hline \vspace{0.3cm}
$(SU(n+1)/T^{n},\cvg)$ & $\lambda_1(t)=1$\\ \vspace{0.3cm}
$(SO(2n+1)/T^{n},\cvh)$ & $\lambda_1(t)=\dfrac{n}{2n-1}$\\ \vspace{0.3cm}
$(Sp(n)/T^{n},\cvk)$ & $\dfrac{4n-1}{4(n+1)}\leq \lambda_1(t)\leq 1$\\ \vspace{0.3cm}
$(SO(2n)/T^{n},\cvm)$ & $\lambda_1(t)= 1$\\ \vspace{0.3cm}
$(G_2/T,\cvn)$ & $\dfrac{1}{2}\leq \lambda_1(t)\leq \dfrac{7}{6}$\\ 
\end{tabular}
\end{table}
\end{proposition}

{\bf Proof:} By Corollaries \ref{primeiroaut1} and \ref{primeiroaut2} in Section \ref{Section2.2}, $$\mu_1\leq \lambda_1(t)\leq \beta_1,$$for all $0< t\leq 1$, where $\lambda_1(t),\mu_1$ and $\beta_1$ are the first positive eigenvalues of $\Delta_t, \Delta_g$ and $\Delta_h$, which are the Laplacians on the total space of the canonical variation, on the original total space and on the basis, respectively. From the above description of the spectra of the Laplacians on maximal flag manifolds and isotropy irreducible Hermitian symmetric spaces, we have all the first positive eigenvalues of the Laplacians both on the original total spaces and on the base, i.e, we have the values of $\beta_1$ and $\mu_1$, which allow us applying the inequality $\mu_1\leq \lambda_1(t)\leq \beta_1$ in order to obtain the respective values and estimates for $\lambda_1(t)$ given in the table above.  \cqd 
\section{Bifurcation and Local Rigidity Instants for
Canonical Variations on Maximal Flag
Manifolds}
\label{section5}

We will now prove the main results of this work. We will determine the bifurcation and local rigidity instants for the canonical variations $\cvg,\cvh,\cvk, \cvm, \cvn$ defined in the Section \ref{Section2.3}. The criterion used to find such instants is based on comparison between an expression of the eigenvalues of the Laplacian relative to the respective canonical variation and a multiple of its scalar curvature. This method is related to the notion of Morse index.   

We will obtain an expression for the scalar curvature of each canonical variation described above according to a general formula due M. Wang and W. Ziller in \cite{WZ} that can be applied to calculate the scalar curvature of homogeneous metrics on reductive homogeneous spaces. %From this formula, we can see clearly that these kind of metrics have constant scalar curvature and are trivial solutions to the Yamabe problem.

Thereafter, necessary conditions for classifications of bifurcation and local rigidity instants in the interval $]0,1[$ can be deduced by using the expressions of the scalar curvature in all cases.

\subsection{Scalar Curvature}
The general formula of the scalar curvature for reductive  homogeneous spaces which we use here is obtained in \cite{WZ}.

Let $G$ be a compact connected Lie group and $K\subset G$ a closed subgroup of $G$. Assume that $K$ is connected, which corresponds to the case that $G/K$ is simply connected. Let $\m$ be the $(-B)$-orthogonal complement to $\sub$ in $\lie$, where $\lie$ and $\sub$ are the Lie algebras of $G$ and $K$, respectively, and $B$ is the Cartan-Killing form of $\lie$. It is known that the isotropy representation $\m$ of $K$ decomposes into a direct sum of inequivalent irreducible submodules, $$\m=\m_1\oplus\m_2\oplus\ldots\oplus\m_r,$$ with $\m_1,\m_2,\ldots,\m_r$ such that $\Ad(K)\m_i\subset\m_i$, for all $i=1,\ldots,r$. Thus, each $G$-invariant metric on $G/K$ can be represented by a inner product on $\m$ given by  $$t_1(-B)|_{\m_1}+t_2(-B)|_{\m_2}+\ldots+t_r(-B)|_{\m_r}, \ \ t_i>0, \ \ i=1,\ldots,r.$$ 

Let ${X_{\alpha}}$ be a $(-B)$-orthonormal basis adapted to the $\Ad(K)$-invariant decomposition of $\m$, i.e., $X_{\alpha}\in \m_i$ for some $i$, and $\alpha<\beta$ if $i<j$ with $X_{\alpha}\in\m_i$ and $X_{\beta}\in\m_j$. Set $A^{\gamma}_{\alpha\beta}=-B([X_{\alpha},X_{\beta}],X_{\gamma})$, so that $[X_{\alpha},X_{\beta}]_{\m}=\displaystyle\sum_{\gamma}{A^{\gamma}_{\alpha\beta}}X_{\gamma},$ and define 

$$\left[
\begin{array}{cc}
k\\
ij\\
\end{array}
\right]=\sum({A^{\gamma}_{\alpha\beta}})^2,$$ where the sum is taken over all indices $\alpha,\beta,\gamma$ with $X_{\alpha}\in\m_i,X_{\beta}\in\m_j $ and $X_{\gamma}\in\m_k$. Note that $\left[
\begin{array}{cc}
k\\
ij\\
\end{array}
\right]$ is independent of the $(-B)$-orthogonal basis chosen for $\m_i,\m_j,\m_k$, but it depends on the choice of the decomposition of $\m$. In addition, $\left[
\begin{array}{cc}
k\\
ij\\
\end{array}
\right]$ is continuous function on the space of all $(-B)$-orthogonal ordered decomposition of $\m$ into $\Ad(K)-$irreducible summands and also is non-negative and symmetric in all 3 indices. The set $\{X_{\alpha}/\sqrt{t_i};X_{\alpha}\in\m_i\}$ is a orthonormal basis of $\m$ with respect to $$\left\langle \cdot,\cdot\right\rangle=t_1(-B)|_{\m_1}+t_2(-B)|_{\m_2}+\ldots+t_r(-B)|_{\m_r}.$$Then the scalar curvature of $g$ determined by $\left\langle \cdot,\cdot\right\rangle$ is
\begin{equation}
\textrm{scal}(g)=\frac{1}{2}\displaystyle\sum_{l=1}^{r}{\frac{d_l}{t_l}-\frac{1}{4}\displaystyle\sum_{i,j,k}\left[
\begin{array}{cc}
k\\
ij\\
\end{array}
\right]\frac{t_k}{t_it_j}}, \label{curvesc}
\end{equation}
for all $x\in G/K$, $d_i=\dim\m_i$, $i=1,\ldots,r$. See \cite{WZ} for details.

We will use now the above formula in order to compute the scalar curvature $\text{scal}(t)$ of each 1-parameter family of homogeneous metrics $\cvg,\cvh,\cvk, \cvm, \cvn$ defined on some of the maximal flag manifold associated with one of the classical simple Lie groups $SU(n+1),SO(2n+1), Sp(n), SO(2n)$ and $G_2$, respectively.

First, are necessary some observations about the numbers $\left[
\begin{array}{cc}
k\\
ij\\
\end{array}
\right]$. For a maximal flag manifold $G/K$, we have the $(-B)$-orthogonal decomposition $\displaystyle\sum_{\alpha\in R^+}\m_{\alpha}$ of its isotropy representation $\m$, where $B$ is the Cartan-Killing form of $\lie$ and $\m_{\alpha}=\R A_{\alpha}+\R S_{\alpha}$, with the vectors $ A_{\alpha}$ and $S_{\alpha}$ defined by the elements of a Weyl basis of $\lie=\text{Lie}(G)$. This allows us rewrite the above splitting of $\m$ as $\m=\m_1\oplus\ldots\oplus\m_s$, where $s=\left|R^+\right|$. Since $B( A_{\alpha},A_{\alpha})=B( S_{\alpha},S_{\alpha})=-2$ and $B( A_{\alpha},S_{\alpha})=0$, the set 
\begin{equation}
\left\{X_{\alpha}=\dfrac{A_{\alpha}}{\sqrt{2}},Y_{\alpha}=\dfrac{S_{\alpha}}{\sqrt{2}};\alpha\in R^+\right\}, \label{BW} 
\end{equation}
is a $(-B)$-orthonormal basis of $\m$. If we denote for simplicity such a basis by $e_{\alpha}=\{A_{\alpha},S_{\alpha}\}$, then the notation $\left[
\begin{array}{cc}
k\\
ij\\
\end{array}
\right]$ can be rewritten as $\left[
\begin{array}{cc}
\gamma\\
\alpha\beta\\
\end{array}
\right]$ where $e_{\alpha}$, $e_{\beta}$ and $e_{\gamma}$ are the $(-B)$-orthogonal bases of the modules $\m_{\alpha}$, $\m_{\beta}$ and $\m{\gamma}$, respectively.

\begin{remark}\label{remarkscal} Recall that, if $\alpha,\beta\in R$ such that $\alpha+\beta\neq 0$, then root spaces satisfy $[\lie_{\alpha},\lie_{\beta}]=\lie_{\alpha+\beta}$ and $B(\lie_{\alpha},\lie_{\beta})=0.$ Since $\left[
\begin{array}{cc}
\gamma\\
\alpha\beta\\
\end{array}
\right]\neq 0$ if and only if $B([\m_{\alpha},\m_{\beta}],\m{\gamma})\neq 0$, we can conclude that $\left[
\begin{array}{cc}
\gamma\\
\alpha\beta\\
\end{array}
\right]\neq 0$ if and only if the positive roots $\alpha,\beta$ and $\gamma$ are such that $\alpha+\beta-\gamma=0$. Thus, we have $\left[
\begin{array}{cc}
\alpha+\beta\\
\alpha\beta\\
\end{array}
\right]\neq 0$, for any $\alpha,\beta$ such that $\alpha+\beta\in R$. Moreover, for $\left[
\begin{array}{cc}
\gamma\\
\alpha\beta\\
\end{array}
\right]$, we will have in the second summation of \eqref{curvesc} the following: 
\begin{enumerate}
\item[(a)] $\left[
\begin{array}{cc}
\gamma\\
\alpha\beta\\
\end{array}
\right]$ $\times $ $\dfrac{1}{t^2}$ when $t^2(-B)|_{\m_{\alpha}},t^2(-B)|_{\m_{\beta}},t^2(-B)|_{\m_{\gamma}}$ are vertical components and $\alpha+\beta$ is a positive root, 
\item[(b)]$\left[
\begin{array}{cc}
\gamma\\
\alpha\beta\\
\end{array}
\right]$ $\times$ $\dfrac{1}{t^2}$ when $t^2(-B)|_{\m_{\alpha}}$ is vertical, $(-B)|_{\m_{\beta}},(-B)|_{\m_{\gamma}}$ are horizontal components and $\alpha+\beta$ is a positive root,
\item[(c)]$\left[
\begin{array}{cc}
\alpha\\
\gamma\beta\\
\end{array}
\right]$ $\times$ $t^2$ when $(-B)|_{\m_{\alpha}}$ is vertical, $(-B)|_{\m_{\beta}},(-B)|_{\m_{\gamma}}$ are horizontal components and $\alpha+\beta$ is a positive root
\item[(d)] $\left[
\begin{array}{cc}
\gamma\\
\alpha\beta\\
\end{array}
\right]=0$ $\times$ $t^2$ otherwise.
\end{enumerate}
\end{remark}

%\begin{proposition}[\cite{Arvanit1}] For a maximal flag manifold $G/K$ the triples $\left[
%\begin{array}{cc}
%\alpha+\beta\\
%\alpha\beta\\
%\end{array}
%\right]\neq 0$ are given by $$\left[
%\begin{array}{cc}
%\alpha+\beta\\
%\alpha\beta\\
%\end{array}
%\right]= 2m_{\alpha,\beta}^2,$$ where $m_{\alpha,\beta}^2$ are the structure constants of the Weyl basis of $\lie=\emph{Lie}(G)$.
%\end{proposition}

%\begin{remark}[\cite{Arvanit1}]If $\alpha,\beta\in R$ are such that $\alpha-\beta \in R$, similarly one obtains $$\left[
%\begin{array}{cc}
%\alpha-\beta\\
%\alpha\beta\\
%\end{array}
%\right]= 2m_{\alpha,-\beta}^2.$$
%\end{remark}

We will now compute the scalar curvature of $\cvg$, defined on $SU(n+1)/T^{n}$.

\begin{lemma}[\cite{sakane}]\label{simbsu} For $SU(n+1)/T^{n}$, considering the decomposition \eqref{decompsu} of the isotropy representation $\m$ of $T^{n}$, one has 
\begin{equation} 
\left[
\begin{array}{cc}
\alpha+\beta\\
\alpha\beta\\
\end{array}
\right]=\left\{
\begin{array}{rcl}
\dfrac{1}{n+1},& \mbox{if}&k\neq i,j\\ 
0, &\mbox{otherwise}
\end{array}
\right.,
\end{equation}where $\alpha=\lambda_i-\lambda_k$ and $\beta=\lambda_k-\lambda_j$ are positive roots of the root system $$R=\{\alpha_{ij}=\pm(\lambda_i-\lambda_j); i\neq j\}$$ of the Cartan subalgebra $\sub^{\mathbb{C}}$ relative to $\mathfrak{su}(n+1)^{\mathbb{C}}$, $\lambda_i$ given by $\lambda_i(\text{diag}(a_1,\ldots,a_{n+1}))=a_i$, for each $1\leq i\leq n+1$ and $$\sub=\left\{\sqrt{-1}\cdot\text{diag}(a_1,\ldots,a_{n+1});\displaystyle\sum_{i=1}^n{a_i}=0\right\}.$$
\end{lemma}
%%%%%%%%%%%%%%%%%%%%%%%%%%%%%%%%%%%%%%%%%%%%%%%%%%%%%%%Curvatura escalar SU(n+1)%%%%%%%%%%%%%%%%%%%%%%%%%%%%%%%%%%%%%%%%%%%%%%%%%%%
\begin{proposition} \emph{[$\bf {(SU(n+1)/T^{n},\cvg), n\geq 2}$]} Considering the hypothesis of the Lemma \ref{simbsu}, let \linebreak$(SU(n+1)/T^{n},\cvg)$ be the canonical variation of $(SU(n+1)/T^{n},g)$, where $g$ is the normal metric on $SU(n+1)/T^{n}$. Then, the function $\emph{scal}(t)$, for each $t>0$, 
\begin{equation}
\emph{scal}(t)=\dfrac{-2n+n^2(n+1)+4n(n+1) t^2+n(1-n)t^4}{4 (n+1) t^2}\label{scalsu}
\end{equation}
gives the scalar curvature of $\cvg$. 
\end{proposition}

{\bf Proof:} The 1-parameter family $\cvg$ is obtained by multiplying the components $$(-B)|_{\m_{\alpha{ij}}}, 1\leq i<j\leq n,$$ by $t^2$, that is, by scaling in the direction of the fibers by $t^2$ the normal metric $g$ given by the inner product $$g_{eT^{n}}=\displaystyle\sum_{i<j\leq n+1}(-B)|_{\m_{\alpha_{ij}}}.$$  

It is known that one has $\dfrac{n(n-1)}{2}$ vertical components, i.e, $$\left|\{(-B)|_{\m_{\alpha{ij}}}, 1\leq i<j\leq n\}\right|\\=\dfrac{n(n-1)}{2},$$ and $\left|\{(-B)|_{\m_{\alpha{ij}}}, 1\leq i<j= n+1\}\right|=n$ horizontal components of the normal metric $g.$ In addition, the dimension of each submodule $\m_{\alpha}$ is equal to 2. Therefore, the first sum in \eqref{curvesc} is $$\frac{1}{2}\displaystyle\sum_{\alpha\in R^+}\dfrac{d_{\alpha}}{t_{\alpha}}=\dfrac{(n-1) n}{2 t^2}+n.$$In order to obtain the second sum
$\displaystyle\sum_{\alpha,\beta,\gamma\in R^+}
\left[
\begin{array}{cc}
\gamma\\
\alpha\beta\\
\end{array}
\right] \dfrac{t_{\gamma}}{t_{\alpha}t_{\beta}}$ in \eqref{curvesc}, with $t_{\gamma},{t_{\alpha},t_{\beta}}>0$ being the coefficients of $(-B)|_{\m_{\alpha}},(-B)|_{\m_{\alpha}},(-B)|_{\m_{\alpha}}$, respectively, it is sufficient to know the number of triples $\left[
\begin{array}{cc}
\gamma\\
\alpha\beta\\
\end{array}
\right]$ multiplying $t^2$ and the number of these constants multiplying $\dfrac{1}{t^2}$. In fact, $\dfrac{t_{\gamma}}{t_{\alpha}t_{\beta}}=t^2$ or $\dfrac{1}{t^2}$, since the coefficients of the canonical variation $g_t$ are $t_{\alpha}=t^2$ or $t_{\alpha}=1$, depending on whether the component $(-B)|_{\m_{\alpha}}$ is vertical or horizontal. Furthermore, the nonzero triples are equal to the same value, namely $0\neq\left[
\begin{array}{cc}
\gamma\\
\alpha\beta\\
\end{array}
\right]=\dfrac{1}{n+1}$ according Lemma \ref{simbsu}.

Let us now determine the number of triples for each of the cases (a), (b) and (c) in Remark \ref{remarkscal} above. \\
CASE (a): Let $N_1$ the total number of triples in (a). Using the above notations, \linebreak$\left[
\begin{array}{cc}
\gamma\\
\alpha\beta\\
\end{array}
\right]=\left[
\begin{array}{cc}
\lambda_i-\lambda_j\\
\lambda_i-\lambda_k\lambda_k-\lambda_j\\
\end{array}
\right],$ with $\alpha=\lambda_i-\lambda_k,$ $\beta=\lambda_k-\lambda_j$ and $\gamma=\alpha+\beta=\lambda_i-\lambda_j,$ $i<k<j\leq n$, since $(-B)|_{\m_{\alpha}},(-B)|_{\m_{\beta}},(-B)|_{\m_{\gamma}}$ are vertical components. Fixed $i=1$, we will have $\dfrac{(n-1)(n-2)}{2}$ symbols of the type $\left[
\begin{array}{cc}
\lambda_1-\lambda_j\\
\lambda_1-\lambda_k\lambda_k-\lambda_j\\
\end{array}
\right]\neq 0;$ if we fix $i=2$, one has $\dfrac{(n-2)(n-3)}{2}$ symbols of the type $\left[
\begin{array}{cc}
\lambda_2-\lambda_j\\
\lambda_2-\lambda_k\lambda_k-\lambda_j\\
\end{array}
\right]\neq 0,$ and so on, getting $\dfrac{(n-i)(n-i-1)}{2}$ symbols of the type $\left[
\begin{array}{cc}
\lambda_i-\lambda_j\\
\lambda_i-\lambda_k\lambda_k-\lambda_j\\
\end{array}
\right]\neq 0$ for each $1\leq i\leq n-2$. It follows that the number of symbols in (a), not counting the permutations, is $$\displaystyle\sum_{i=1}^{n-2}\dfrac{1}{2}(n-i)(n-i-1)=\dfrac{1}{6}(2 n - 3 n^2 + n^3).$$Recall that in (a), $\left[
\begin{array}{cc}
\lambda_i-\lambda_j\\
\lambda_i-\lambda_k\lambda_k-\lambda_j\\
\end{array}
\right]=\left[
\begin{array}{cc}
\gamma\\
\alpha\beta\\
\end{array}
\right]$,with  $i<k<j\leq n$. Then, $t_{\alpha},t_{\beta},t_{\gamma}$ are equal to $t^2$, which implies that $\dfrac{t_{\gamma}}{t_{\alpha}t_{\beta}}=\dfrac{1}{t^2}$, and by symmetry of $\left[
\begin{array}{cc}
\gamma\\
\alpha\beta\\
\end{array}
\right],$ we must considering six times the number of symbols above to obtain the total number of these symbols in (a), hence $$N_1=6\cdot\dfrac{1}{6}(2 n - 3 n^2 + n^3)=(2 n - 3 n^2 + n^3).$$CASE (b): Let $N_2$ the total number of triples in (b); $\left[
\begin{array}{cc}
\gamma\\
\alpha\beta\\
\end{array}
\right]=\left[
\begin{array}{cc}
\lambda_i-\lambda_j\\
\lambda_i-\lambda_k\lambda_k-\lambda_j\\
\end{array}
\right],$ with $\alpha=\lambda_i-\lambda_k,$ $\beta=\lambda_k-\lambda_j$ and $\gamma=\alpha+\beta=\lambda_i-\lambda_j,$ $i<k<j=n+1$, since $(-B)|_{\m_{\alpha}}$ is vertical and $(-B)|_{\m_{\beta}},(-B)|_{\m_{\gamma}}$ are horizontal components. Fixed $i=1$, we will have $(n-1)$ triples of the type $\left[
\begin{array}{cc}
\lambda_1-\lambda_j\\
\lambda_1-\lambda_k\lambda_k-\lambda_j\\
\end{array}
\right]\neq 0;$ if we fix $i=2$, one has $(n-2)$ triples of the type $\left[
\begin{array}{cc}
\lambda_2-\lambda_j\\
\lambda_2-\lambda_k\lambda_k-\lambda_j\\
\end{array}
\right]\neq 0,$ and so on, getting $(n-i)$ triples of the type $\left[
\begin{array}{cc}
\lambda_i-\lambda_j\\
\lambda_i-\lambda_k\lambda_k-\lambda_j\\
\end{array}
\right]\neq 0$ for each $1\leq i\leq n-1$. It follows that the number of triples in (b), not counting the permutations, is $$\displaystyle\sum_{i=1}^{n-1}i=\dfrac{n(n-1)}{2}.$$In (b), $\left[
\begin{array}{cc}
\lambda_i-\lambda_j\\
\lambda_i-\lambda_k\lambda_k-\lambda_j\\
\end{array}
\right]=\left[
\begin{array}{cc}
\gamma\\
\alpha\beta\\
\end{array}
\right]$ then, $t_{\alpha}=t^2,t_{\beta}=t_{\gamma}=1$, which implies $\dfrac{t_{\gamma}}{t_{\alpha}t_{\beta}}=\dfrac{t_{\gamma}}{t_{\beta}t_{\alpha}}=\dfrac{t_{\beta}}{t_{\gamma}t_{\alpha}}=\dfrac{t_{\beta}}{t_{\alpha}t_{\gamma}}=\dfrac{1}{t^2}$. By symmetry of $\left[
\begin{array}{cc}
\gamma\\
\alpha\beta\\
\end{array}
\right],$ we must considering four times the number of triples above to obtain the total number $N_2$ of these triples in (b), hence  $$N_2=4\cdot\dfrac{n(n-1)}{2}=2n(n-1).$$    
CASE (c): Similarly to the previous case, not counting permutations, we have $\dfrac{n(n-1)}{2}$ symbols $\left[
\begin{array}{cc}
\alpha\\
\gamma\beta\\
\end{array}
\right]$ multiplying $\dfrac{t_{\alpha}}{t_{\gamma}t_{\beta}}=\dfrac{t_{\alpha}}{t_{\beta}t_{\gamma}}=t^2$. By symmetry, in this case, we must consider two times the number above in order to obtain the total number $N_3$ of triples $\left[
\begin{array}{cc}
\alpha\\
\gamma\beta\\
\end{array}
\right],$ which multiply in (c), i.e, $N_3=2\cdot \dfrac{n(n-1)}{2}=n(n-1)$.

Therefore, since $\left[
\begin{array}{cc}
\alpha\\
\gamma\beta\\
\end{array}
\right]=\dfrac{1}{n+1}$ whenever $\left[
\begin{array}{cc}
\alpha\\
\gamma\beta\\
\end{array}
\right]\neq 0,$  
\begin{eqnarray*}
\text{scal}(t)&=&\frac{1}{2}\displaystyle\sum_{\alpha\in R^+}\dfrac{d_{\alpha}}{t_{\alpha}}-\dfrac{1}{4}\displaystyle\sum_{\alpha,\beta,\gamma\in R^+}
\left[
\begin{array}{cc}
\gamma\\
\alpha\beta\\
\end{array}
\right] \dfrac{t_{\gamma}}{t_{\alpha}t_{\beta}}\\
             &=&\dfrac{(n-1) n}{2 t^2}+n-\dfrac{1}{4}\left(\dfrac{N_1}{(n+1)t^2}+\dfrac{N_2}{t^2(n+1)}+\dfrac{N_3t^2}{n+1}\right)\\
						 &=&\dfrac{-2n+n^2(n+1)+4n(n+1) t^2+n(1-n)t^4}{4 (n+1) t^2}.
\end{eqnarray*}						
\cqd
%\lipsum[12]
%See awesome Table~\ref{tab:table}.

%\begin{table}
 %\caption{Sample table title}
  %\centering
  %\begin{tabular}{lll}
   % \toprule
   % \multicolumn{2}{c}{Part}                   \\
   % \cmidrule(r){1-2}
    %Name     & Description     & Size ($\mu$m) \\
    %\midrule
    %Dendrite & Input terminal  & $\sim$100     \\
    %Axon     & Output terminal & $\sim$10      \\
    %Soma     & Cell body       & up to $10^6$  \\
    %\bottomrule
  %\end{tabular}
  %\label{tab:table}
%\end{table}

The procedure to obtain the formulae for $\text{scal}(t)$ of the other canonical variations is analogous to the previous case and we omit their proof. 

\begin{proposition} \emph{[$\bf {(SO(2n+1)/T^n,\cvh), n\geq 4}$]} Let $SO(2n+1)/T^n,\cvh)$ be the canonical variation of $(SO(2n+1)/T^n,g)$, where $g$ is the normal metric on $SO(2n+1)/T^n$. Then, the function $\emph{scal}(t),$ that for each $t>0$ gives the scalar curvature of $\cvh,$ is given by 
\begin{equation}
\emph{scal}(t)=\dfrac{5 n^3-2 n^2 t^4+8 n^2 t^2-7 n^2+2 n t^4-4 n t^2+2 n}{4 (2 n-1) t^2}.\label{scalso1}
\end{equation}
\end{proposition}

%The next case is the scalar curvature $\text{scal}(t)$ of the canonical variation of the maximal flag manifold associated with the complex simple Lie algebra $\mathfrak{sp}(n)^{\mathbb{C}}$, endowed with the normal metric.

\begin{proposition} \emph{[$\bf {(Sp(n)/T^n,\cvk), n\geq 3}$]} Let $(Sp(n)/T^n,\cvk)$ be the canonical variation of $(Sp(n)/T^n,g)$, where $g$ is the normal metric on $Sp(n)/T^n$. Then, the function $\emph{scal}(t),$ that for each $t>0$ gives the scalar curvature of $\cvk,$ is given by 
\begin{equation}
\emph{scal}(t)=\dfrac{-2 n^3 t^4+24 n^3 t^2+5 n^3+48 n^2 t^2+9 n^2+2 n t^4+24 n t^2-14 n}{24 (n+1) t^2}.\label{scalsp}
\end{equation}
\end{proposition}

%The last expression for $\text{scal}(t)$ of the canonical variations for the maximal flag manifolds associated with a complex classical simple Lie algebra is the one associated with $\mathfrak{so}(2n)^{\mathbb{C}}$ provided with its normal metric.

\begin{proposition} \emph{[$\bf {(SO(2n)/T^n,\cvm), n\geq 4}$]} Let $(SO(2n)/T^n,\cvm)$ be the canonical variation of $(SO(2n)/T^n,g)$, where $g$ is the normal metric on $SO(2n)/T^n$. Then, the function $\emph{scal}(t),$ that for each $t>0$ gives the scalar curvature of $\cvm,$ is given by 
\begin{equation}
\emph{scal}(t)=\dfrac{-2 n^2 t^4+24 n^2 t^2+5 n^2+4 n t^4-24 n t^2+2 n}{24 t^2}.\label{scalso2}
\end{equation}
\end{proposition}

%Consider now the formula for $\text{scal}(t)$ of the canonical variation $\cvn$ of the normal metric on the maximal flag manifold associated with the exceptional simple Lie algebra of the type $G_2$. %Before, we will enunciate the following useful lemma.
%
%\begin{lemma}[\cite{Arvanit1}]\label{lemaG2} For the isotropy representation $\m=\m_{\alpha_1}\oplus\m_{\alpha_2}\oplus\m_{\alpha_1+\alpha_2}\oplus\m_{\alpha_1+2\alpha_2}\oplus\m_{\alpha_1+3\alpha_2}\oplus\m_{2\alpha_1+3\alpha_2}$ of $G_2/T$, put $\m_1=\m_{\alpha_1},$ $\m_2=\m_{\alpha_2},$ $\m_3=\m_{\alpha_1+\alpha_2},$ $\m_4=\m_{\alpha_1+2\alpha_2},$ $\m_5=\m_{\alpha_1+3\alpha_2}$ and $\m_6=\m_{2\alpha_1+3\alpha_2}$. Using the current notation, the non zero triples $\left[
%\begin{array}{cc}
%k\\
%ij\\
%\end{array}
%\right] $ of the maximal flag $G_2/T$ are $$\left[
%\begin{array}{cc}
%3\\
%12\\
%\end{array}
%\right]=\left[
%\begin{array}{cc}
%5\\
%24\\
%\end{array}
%\right]=\left[
%\begin{array}{cc}
%6\\
%34\\
%\end{array}
%\right]=\left[
%\begin{array}{cc}
%6\\
%15\\
%\end{array}
%\right]=\dfrac{1}{4}, \ \ \emph{and} \ \ \left[
%\begin{array}{cc}
%4\\
%23\\
%\end{array}
%\right]=\dfrac{1}{3}.$$
%\end{lemma}

\begin{proposition}\label{propscalG} \emph{[$\bf {(G_2/T,\cvn)}$]} Let $(G_2/T,\cvn)$ be the canonical variation of $(G_2/T,g)$, where $g$ is the normal metric on $G_2/T$. Then, the function $\emph{scal}(t),$ that for each $t>0$ gives the scalar curvature of $\cvn,$ is given by 
\begin{equation}
\emph{scal}(t)=\dfrac{2+12t^2-2t^4}{3t^2}.\label{scalG}
\end{equation}
\end{proposition}
%{\bf Proof:} Straightforward computation from the formula \eqref{curvesc}.
\subsection{Bifurcation and Local Rigidity}

If the Morse index changes when passing a degeneracy instant this is actually a bifurcation instant. We will apply such principle in order to find bifurcation instants for the canonical variations $\cvg,\cvh,\cvk,\cvm,\cvn, 0<t<1$. 

It was defined previously that a {\it degeneracy instant} $t_{\ast}>0$ for $g_t$ in ${\mathcal{R}}^k(M)$, with $g_1=g$, is an instant such that $\dfrac{\textrm{scal}(g_{t_{\ast}})}{m-1}\in\sigma(\Delta_{t_{\ast}})$.

Denoting by $$\{0<\lambda^g_1\leq \lambda^g_2\leq\ldots\leq\lambda^g_j\leq\ldots\}$$ the sequence of positive eigenvalues of $\Delta_g$, the Morse index of a Riemannian metric $g$ is $$N(g)=\max\left\{j\in\mathbb{N};\lambda^g_j<\dfrac{\textrm{scal}(g)}{m-1}\right\},$$ where is the Laplacian $\Delta_g$ acting on $C^{\infty}(M)$, $M$ provide with the Riemannian metric $g$ and $m=\dim M.$ 

The following result is a sufficient condition for a degeneracy instant $0<t_{\ast}<1$ to be a bifurcation instant when $\dfrac{\text{scal}(t_{\ast})}{m-1}$ is a constant eigenvalue of the Laplacian $\Delta_{t_{\ast}}$. 

\begin{proposition}[\cite{pacificjournal}]\label{propbifurcation} Let $(M,g)$ be a closed Riemannian manifold with $\dim M\geq 3$ and $\pi:(M,g)\longrightarrow (B,h)$ a Riemannian submersion with totally geodesic fibers isometric to $(F,\kappa)$, where $\dim F\geq 2$ and $\emph{scal}(F)>0$. Denote by $\lambda\in\sigma(\Delta_h)\subset\sigma(\Delta_g)$ a constant eigenvalue of $\Delta_{t_{\ast}}$ such that $\dfrac{\emph{scal}(g_{t_{\ast}})}{m-1}=\lambda$, $g_{t_{\ast}}$ canonical variation of $g$ at $0<t_{\ast}<1$ and $\Delta_{t_{\ast}}$ the Laplacian on $(M,g_{t_{\ast}})$. If $$\dfrac{\emph{scal}(g_{t_{\ast}})}{m-1}<\lambda^{1,1}(t_{\ast})=\mu_1+(\frac{1}{t_{\ast}^2}-1)\phi_1,$$ $\mu_1\in \sigma(\Delta_g)$ the first positive eigenvalue of $\Delta_g$ and $\phi_1\in \sigma(\Delta_v)$ the first positive eigenvalue of the vertical Laplacian $\Delta_v$, then $t_{\ast}$ is a bifurcation instant for $g_t$.
\end{proposition}

{\bf Proof:} It is sufficient to show that the Morse index changes when passing the degeneracy value $0<t_{\ast}<1$, i.e, for $\epsilon>0$ sufficiently small, $N(g_{t_{\ast}-\epsilon})\neq N(g_{t_{\ast}+\epsilon})$ then $t_{\ast}$ is a bifurcation instant. 

Observe that, if the Morse index does not change, there must be a compensation of eigenvalues. Namely, there must exist nonconstant eigenvalues $\lambda^{k_1,j_1}(t),\ldots,\lambda^{k_n,j_n}(t)$ of $\Delta_t$, whose combined multiplicity equals the multiplicity
of $\lambda$, such that, $$\lambda<\dfrac{\textrm{scal}(g_t)}{m-1}<\lambda^{k_i,j_i}(t), \ \ \forall t<t_{\ast} \ \ \textrm{(close to $t_{\ast}$) and}\ \ 1\leq i\leq n,$$ 
$$\lambda>\dfrac{\textrm{scal}(g_t)}{m-1}>\lambda^{k_i,j_i}(t), \ \ \forall t>t_{\ast} \ \ \textrm{(close to $t_{\ast}$) and}\ \ 1\leq i\leq n.$$
If $\dfrac{\textrm{scal}(g_{t_{\ast}})}{m-1}<\mu_1+(\dfrac{1}{t_{\ast}^2}-1)\phi_1,$ $\mu_1\in \sigma(\Delta_g)$ first positive eigenvalue of $\Delta_g$ and $\phi_1\in \sigma(\Delta_v)$ first positive eigenvalue of the vertical Laplacian $\Delta_v$, then $$\lambda=\dfrac{\textrm{scal}(g_{t_{\ast}})}{m-1}<\mu_1+\left(\frac{1}{t_{\ast}^2}-1\right)\phi_1<\lambda^{k,j}(t_{\ast}), \ \ \forall \mathbb{Z}\ni k,j>0,$$ since $\lambda^{k,j}(t_{\ast})=\mu_k+(\dfrac{1}{t_{\ast}^2}-1)\phi_j$, $\mu_1<\mu_k\in\sigma(\Delta_g)$, $\phi_1<\phi_k\in\sigma(\Delta_v)$  and $(\dfrac{1}{t_{\ast}^2}-1)>0$ when $0<t_{\ast}<1$. It follows that every nonconstant eigenvalue $\lambda^{k,j}(t)$ is strictly greater than $\lambda=\dfrac{\textrm{scal}(g_{t_{\ast}})}{m-1}$ for $t$ sufficiently close to $t_{\ast}$, so that there is no compensation of eigenvalue and the Morse index must change when passing $t_{\ast}$. Since $\dfrac{\textrm{scal}(g_{t_{\ast}})}{m-1}=\lambda\in\sigma(\Delta_{t_{\ast}})$, $t_{\ast}$ is a degeneracy instant for $g_t$, and then $t_{\ast}$ is bifurcation instant for $g_t$. \cqd 

%Later on, $g$ is the normal homogeneous metric on the maximal flag manifold $G/T$, represented by $(-B)|_{\m}$, $B$ Cartan Killing form of the Lie group $G$, $\m=\p\oplus\q$ the isotropy representation of the maximal torus $T\subset G$, $H/T$ represents the fibers provided with the induced homogeneous metric $\hat{g}$ represented by the inner product $(-B)|_{\p}$ and the basis space $G/H$ is provided with the symmetric metric $\breve{g}$ induced by $(-B)|_{\q}$, $\p$ and $\q$ vertical and horizontal distributions, respectively.

For each canonical variation introduced above, defined on the maximal flag manifolds associated with a complex classical simple Lie algebra, we obtain the following properties

\begin{lemma}\label{autoconstante}Every degeneracy instant $t_{\ast}$ for $\cvg,\cvh,\cvk,\cvm$ and $\cvn$  is such that $$\dfrac{\emph{scal}(t_{\ast})}{m-1}\in\sigma(\Delta_{\breve{g}})\subset \sigma(\Delta_{t_{\ast}}),$$ with $\sigma(\Delta_{t_{\ast}})$ being the spectrum of the Laplacians on the total spaces $SU(n+1)/T^{n},\linebreak SO(2n+1)/T^{n},Sp(n)/T^{n},SO(2n)/T^{n}$, $G_2/T$ and $\sigma(\Delta_{\breve{g}})$ being the spectrum on the basis spaces $(G/H,\breve{g})$, respectively, $m=\dim G/T.$ In other words, $\dfrac{\emph{scal}({t_{\ast}})}{m-1}$ is eigenvalue of $\Delta_{t_{\ast}}$ if and only if $\dfrac{\emph{scal}({t_{\ast}})}{m-1}$ is  a constant eigenvalue $\lambda^{k,0}(t)\in\sigma(\Delta_{t_{\ast}})$, for some $1\leq k\in\mathbb{Z}.$

\end{lemma}

{\bf Proof:} It is known that an eigenvalue of $\Delta_t$ can be written as $$\lambda^{k,j}(t)=\mu_k+\left(\frac{1}{t^2}-1\right)\phi_j,$$ for some $\mu_k\in\sigma(\Delta_g)$ and $\phi_j\in\sigma(\Delta_{\hat{g}})$. To proof that $\dfrac{\text{scal}({t_{\ast}})}{m-1}$ is a constant eigenvalue $\lambda^{k,0}(t)\in\sigma(\Delta_{t_{\ast}})$, for some $1\leq k\in\mathbb{Z},$ it is enough to show that there are no $k,j>0$ such that $$\dfrac{\text{scal}({t_{\ast}})}{m-1}=\lambda^{k,j}(t)$$ nor $j>0$ such that $$\dfrac{\textrm{scal}(t)}{m-1}=\lambda^{0,j}(t).$$
In order to prove the inequality $$\dfrac{\text{scal}(t)}{m-1}<\lambda^{1,1}(t)=\mu_1+\left(\dfrac{1}{t^2}-1\right)\phi_1,\forall\, 0<t\leq 1,$$ for each of the above (a)-(e) cases we define the function $f$ given by $$f(t)=\dfrac{\text{scal}(t)}{m-1}-\lambda^{1,1}(t)=\dfrac{\text{scal}(t)}{m-1}-\mu_1-\left(\dfrac{1}{t^2}-1\right)\phi_1,\forall\, 0<t\leq 1,$$and we verify that $f$ is strictly negative for $0<t\leq 1$. Indeed, in the case of the canonical variation $\cvg$ on $SU(n+1)/T^n$ we obtain $$f(t)=\frac{n^3-n^2 t^4+4 n^2 t^2+n^2+n t^4+4 n t^2-2 n}{4 (n+1) (n (n+1)-1) t^2}-\frac{1}{t^2}, $$where $f(1)=\frac{n^3+4 n^2+3 n}{4 (n+1) (n (n+1)-1)}-1<0,$ for all $n\geq 2$ and 
\begin{eqnarray*}
\frac{df}{dt}(t)&=&\frac{-4 n^2 t^3+8 n^2 t+4 n t^3+8 n t}{4 (n+1) (n (n+1)-1) t^2}\\
                & & -\frac{n^3-n^2 t^4+4 n^2 t^2+n^2+n t^4+4 n t^2-2 n}{2 (n+1) (n (n+1)-1)t^3}+\frac{2}{t^3}>0,
\end{eqnarray*}
for $0<t\leq 1$. In the other cases, the proof is analogous.
Furthermore, it can be shown, by applying elementary differential calculus, that there is no $0<t<1$ such that $\dfrac{\textrm{scal}(t)}{m-1}=\lambda^{0,j}(t),$ for all integer $j\geq 1$. We will present the proof of this property for the canonical variation $\cvg$. The proof of the another cases are entirely analogous.

From the description of the spectra of maximal flag manifolds and symmetric spaces, the first positive eigenvalue of the Laplacian $\Delta_t$ acting on functions on the total space $(SU(n+1)/T^{n},\cvg)$ is $\lambda_1(t)=1$, for all $0<t\leq 1$. Since $$\textrm{scal}(t)=\dfrac{-2n+n^2(n+1)+4n(n+1) t^2+n(1-n)t^4}{4 (n+1) t^2}, n\geq 2,$$ for all $t>0$, and $m=\dim SU(n+1)/T^{n}=n(n+1)$, we have 
\begin{eqnarray*}
\dfrac{\textrm{scal}(g_t)}{m-1}&=&\frac{n^3-n^2 t^4+4 n^2 t^2+n^2+n t^4+4 n t^2-2 n}{4 (n+1) (n (n+1)-1)t^2}\\                            & < &\lambda_1(t)=1\\
&\Leftrightarrow & \sqrt{\sqrt{\dfrac{4 n^4+17 n^3+26 n^2+16 n+4}{n^2}}-\dfrac{2 \left(n^2+2 n+1\right)}{n}}<t\leq 1, 
\end{eqnarray*}
Thus, we have that $\cvg$ is locally rigidity at all instant 
$$t_{\ast}\in \ \  ]\sqrt{\sqrt{\dfrac{4 n^4+17 n^3+26 n^2+16 n+4}{n^2}}-\dfrac{2 \left(n^2+2 n+1\right)}{n}},1].$$ Hence, there are no degeneracy instants for $\cvg$ in the interval $]b,1],$ with $$b=\sqrt{\sqrt{\dfrac{4 n^4+17 n^3+26 n^2+16 n+4}{n^2}}-\dfrac{2 \left(n^2+2 n+1\right)}{n}}.$$

The fiber of the canonical fibration $(SU(n+1)/T^{n},\cvg)$ is the maximal flag manifold $(SU(n)/T^{n-1},\hat{g})$, with the induced metric $\hat{g}$, represented by the inner product $(-B)|_{\p}$. The first positive eigenvalue $\phi_1$ of the  Laplacian $\Delta_{\hat{g}}$ of such homogeneous metric is equal to 1, i.e., if $\phi\in\sigma(\Delta_{\hat{g}})$, then $\phi\geq 1$.

Define $\varphi_r:]0,1[\longrightarrow \R$ by $$\varphi_r(t)=\dfrac{\textrm{scal}(g_t)}{m-1}-\left(\dfrac{1}{t^2}-1\right)\cdot r,$$ for some fixed $r\geq 1$. We have that $$\dfrac{d}{dt}(\varphi_r(t))=\frac{n^3 (4 r-1)+n^2 \left(8 r-t^4-1\right)+n \left(t^4+2\right)-4 r}{2 (n+1) \left(n^2+n-1\right) t^3}>0,$$and $$\varphi_r(b)=\frac{n ((C_n-5) r+C_n-4)-2 n^2 (r+1)-2 (r+1)}{(C_n-4) n-2 n^2-2}<0,$$$\forall \, 0<t<1, r\geq 1, n\geq 2,$ where $C_n=\sqrt{4 n^2+\frac{4}{n^2}+17 n+\frac{16}{n}+26}.$ It follows that $\varphi_r(t)$ is negative for all $0<t<b$ and $$\dfrac{\textrm{scal}(g_t)}{m-1}<\left(\dfrac{1}{t^2}-1\right)\cdot r,$$ for any $r\geq 1$, that is, there is no $0<t<1$ such that $\dfrac{\textrm{scal}(t)}{m-1}=\lambda^{0,j}(t)=\left(\dfrac{1}{t^2}-1\right)\cdot\phi_j,$ for all $\phi_j\in\sigma(\Delta_{\hat{g}})$ and we complete the proof for the case of the canonical variation $\cvg$.

As well as for the case of $\cvg$, define for $\cvh, \cvk, \cvm$ and $\cvn$ the function $\varphi_r(t)$ as above with $r\geq 1$, since the first positive eigenvalues of the Laplacians on their respective fibers are also equal to 1. The statement follows by checking that $\dfrac{d}{dt}(\varphi_r(t))>0, \forall\, 0<t<1$ and $\varphi_r(b)<0$, where $0<b<1$ is such that $$\dfrac{\textrm{scal}(g_t)}{m-1}<\lambda_1(t), \forall\, b<t<1,$$with $\lambda_1(t)$ first positive eigenvalues of the Laplacians.
 \cqd
 
%Observe that the proof of the last Lemma does not guarantees that the set of all instants $0<t<1$ such that $\dfrac{\textrm{scal}(g_t)}{m-1}<\lambda_1(t)$ is equal to the interval $]b,1],$ for the canonical variations $\cvk$ and $\cvn.$ In fact, we don't have expressions of $\lambda_1(t)$ in such cases, therefore we can not determine the interval of rigidity instants, only subintervals of them. We used the lower bounds of $\lambda_1(t)$ in order to determine a subset $]b,1]$ of local rigidity instants for the canonical variations $\cvk$ and $\cvn.$ 

%Already for the canonical variations $\cvg$, $\cvh$ and $\cvm$ on $SU(n+1)/T^{n}$, $SO(2n+1)/T^n$ and $SO(2n)/T^n$ we determined all $t\in ]0,1]$ such that $\dfrac{\textrm{scal}(g_t)}{m-1}<\lambda_1(t)$, since $\lambda_1(t)$ is constant on $]0,1]$ in these cases. 

The degeneracy instants for the conical variation $\cvg$, $0<t<1$, on $SU(n+1)/T^{n}$ are given in our following theorem.

\begin{theorem} Let $\cvg$ be the above canonical variation on $SU(n+1)/T^{n}$ and take $$b=\sqrt{\sqrt{\dfrac{4 n^4+17 n^3+26 n^2+16 n+4}{n^2}}-\dfrac{2 \left(n^2+2 n+1\right)}{n}}.$$ Thus, the degeneracy instants for $\cvg$ in $]0,1[$ form a decreasing sequence $\{t^{\text{\bf g}}_{q}\}\subset \, ]0,b]$ such that $t^{\text{\bf g}}_{q}\rightarrow 0$ as $q\rightarrow 0$, with $t^{\text{\bf g}}_1=b$ and for $q>1$, 
\begin{equation}
t^{\text{\bf g}}_{q}=\sqrt{\sqrt{f(q)}-g(q)}. \label{seqsu}
\end{equation}
where 
\begin{eqnarray*}
f(q)&=&\frac{1}{n^2(n-1)^2}\cdot(4 n^6 q^2+n^5 \left(8   q^3+8 q^2-8 q+1\right)+\\
    & & 4 n^4 \left(q^4+4 q^3-3 q^2-4 q+1\right)+n^3 \left(8 q^4-8 q^3-24q^2+5\right)+\\
    & & n^2 \left(-4 q^4-16 q^3+4 q^2+8 q+6\right)+8 n q^2 \left(-q^2+q+1\right)+4 q^4)
 \end{eqnarray*}  
and $g(q)=\frac{2 \left(n^3 q+n^2
   \left(q^2+q-1\right)+n \left(q^2-q-1\right)-q^2\right)}{(n-1) n}.$
\end{theorem}

{\bf Proof:} We must determine all $t\in \ \ ]0,b]$ such that $\dfrac{\textrm{scal}(t)}{m-1}\in\sigma(\Delta_t)$, since by the previous theorem $g_t$ is locally rigidity for $b<t<1$. From Lemma \ref{autoconstante}, if $\dfrac{\textrm{scal}(t)}{m-1}$ is a eigenvalue of $\Delta_t$, then $\dfrac{\textrm{scal}(t)}{m-1}\in\Delta_{\breve{g}}$.  

Thus, it remains verify for which instants $0<t<b<1$ one has $\dfrac{\textrm{scal}(t)}{m-1}=\lambda^{k,0}(t)$.  

We have that $\lambda^{k,0}(t)$ is eigenvalue of $\Delta_t$ if and only if $\lambda^{k,0}(t)$ belongs to the spectrum of the Laplacian on the basis $\mathbb{CP}^n=SU(n+1)/S(U(1)\times U(n))$, provide with the symmetric metric $\breve{g}$ represented by the inner product $(\cdot,\cdot)=-B|_{\mathfrak{q}}$, $\mathfrak{q}$ horizontal space of the original fibration. The spectrum of the Laplacian $\Delta_{\breve{g}}$ on $(\mathbb{CP}^n,\breve{g})$ is $$ \sigma(\Delta_{\breve{g}})=\{\beta_q=\dfrac{q(q+n)}{n+1};q\in\mathbb{N}\}.$$
 
The degeneracy instants for $\cvg$ in $]0,b]$ are the real values $t^g_q$, solutions of the equation $$\dfrac{\textrm{scal}(t)}{m-1}=\beta_q=\dfrac{q(q+n)}{n+1}. $$ The explicit formula for $t^{\text{\bf g}}_{q}$, which represents the solution of the above equation in $t$, is presented in \eqref{seqsu}. Note that the constant eigenvalues $\beta_q=\dfrac{q(q+n)}{n+1}$ tend to $+\infty$ as $q\rightarrow \infty$ and, since $\dfrac{\text{scal}(t)}{m-1}$ is continuous and tends to $+\infty$ as $t\rightarrow 0$, $t^{\text{\bf g}}_{q}\rightarrow 0.$ \cqd 
\begin{figure}[h]
\centering
\includegraphics[height=7cm,width=10cm]{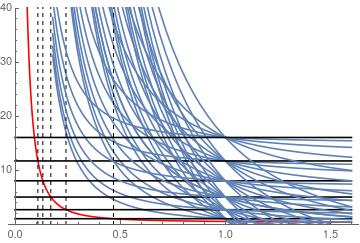}
\caption{For n=2, we have the flag manifold $SU(3)/T^2$, with $m=\dim SU(3)/T^2=6$. Graphics of the function $\dfrac{\textrm{scal}(t)}{m-1}$ in red and of the eigenvalues $\lambda^{k,j}(t)$; the constant eigenvalues in black correspond to $(k,0), 1\leq k\leq 6$, and the non-constant eigenvalues in blue correspond  to $(k,j)$ where $1\leq j\leq k\leq 6$. The dashed vertical lines mark the first five degeneracy instants (which are all bifurcation instants) starting at $b=t^{\text{\bf g}}_1$.}
\end{figure}

%Veremos a seguir que, sob certas hipóteses, é possível garantir a existência de infinitos instantes de degenerescência para uma família $g_t$ obtida a partir de uma submersão Riemanniana $\pi:(M,g)\longrightarrow (B,h)$ com fibras totalmente geodésicas isométricas a $(F,\kappa)$. 

We also compute the Morse index of each $\cvg$ as a critical point of the Hilbert-Einstein functional.

\begin{proposition}\label{morseindexsu}
The Morse index of $\cvg$ is given by 
\begin{equation} 
N(t)=\left\{
\begin{array}{rcl}
\displaystyle\Sigma_{q=1}^r m_q(\mathbb{CP}^n),& \mbox{if}&t^{\text{\bf g}}_{r+1}\leq t < t^{\text{\bf g}}_{r}\\ 
0, &\mbox{if}& t^{\text{\bf g}}_1\leq t\leq 1
\end{array}
\right.,
\end{equation}where $m_q(\mathbb{CP}^n)$ is the multiplicity of the $q$th eigenvalue of the basis $\mathbb{CP}^n=SU(n+1)/S(U(1)\times U(n)).$ 
\end{proposition}

{\bf Proof:} We established in Lemma \ref{autoconstante} that $\dfrac{\textrm{scal}(t)}{m-1}\leq\lambda_1(t)$ for all $t\in [b,1]$, with $b=t^{\text{\bf g}}_1,$ so that there are no eigenvalues of $\Delta_t$ that are less than $\dfrac{\textrm{scal}(t)}{m-1}$ in the interval $[b,1]$. Hence, $N(t)=0$ for $t\in [b,1]$. When $t\rightarrow 0$, whenever $t$ crosses a degeneracy instant $t^{\text{\bf g}}_{q}$, the constant eigenvalue $\lambda^{q,0}(t)$ becomes smaller than $\dfrac{\textrm{scal}(t)}{m-1}$. 

Therefore, the Morse index increases by the multiplicity of $\lambda^{q,0}(t),$ which is the dimension of the corresponding  eigenspace $E^{0}_{q}$. This dimension is given by the Weyl dimension formula and is positive. Thus, $\dim E^{0}_{q}$ also coincides with the multiplicity of the complex projective space $\mathbb{CP}^n,$ concluding the proof. \cqd

The degeneracy instants for the canonical variation $\cvh$, $0<t\leq 1$, on $SO(2n+1)/T^n$ are given in our following theorem.

\begin{theorem} Let $\cvh$ be the above canonical variation on $SO(2n+1)/T^n$ and take $$b=\sqrt{\dfrac{\sqrt{8 n^2+5 n-2}}{\sqrt{2}}-2 n}.$$ Thus, the degeneracy instants for $\cvh$ in $]0,b]$ form a decreasing sequence $\{t^{\text{\bf h}}_{q}\}\subset \, ]0,b]$ such that $t^{\text{\bf h}}_{q}\rightarrow 0$ as $q\rightarrow 0$, with $t^{\text{\bf h}}_1=b$ and for $q>1$, 
\begin{equation}
t^{\text{\bf h}}_{q}=\sqrt{\sqrt{f(q)}+g(q)}. \label{seqso1}
\end{equation}
where 
\begin{eqnarray*}
f(q)&=&\frac{1}{(n-1)^2 n^2}(10 n^5-8 n^4+2 n^3+\left(4 n^4-4      n^2+1\right) q^4 \\ 
    & & +\left(16 n^5-8 n^4-16 n^3+8 n^2+4 n-2\right) q^3 \\  
    & & +\left(-32 n^5+32 n^4+8 n^3-16 n^2+4
   n\right) q \\ 
    & & +\left(16 n^6-16 n^5-28 n^4+24 n^3+8 n^2-8 n+1\right) q^2)
\end{eqnarray*}  
and $g(q)=\dfrac{-4 n^3 q-2 n^2 q^2+2 n^2 q+4 n^2+2 n q-2 n+q^2-q}{2 (n-1) n}.$
\end{theorem}

{\bf Proof:} We must determine all $t\in \ \ ]0,b]$ such that $\dfrac{\textrm{scal}(t)}{m-1}\in\sigma(\Delta_t)$, since by the previous theorem $\cvh$ is locally rigidity for $b<t<1$. From Lemma \ref{autoconstante}, if $\dfrac{\textrm{scal}(t)}{m-1}$ is a eigenvalue of $\Delta_t$, then $\dfrac{\textrm{scal}(t)}{m-1}\in\Delta_{\breve{g}}$.

Thus, it remains verify for which instants $0<t<b<1$ one has $\dfrac{\textrm{scal}(t)}{m-1}=\lambda^{k,0}(t)$. 

From Proposition \ref{constanteigen}, Section \ref{Section2.2}, $\lambda^{k,0}(t)$ is eigenvalue of $\Delta_t$ if and only if $\lambda^{k,0}(t)$ belongs to the spectrum of the Laplacian on the basis $S^{2n}=SO(2n+1)/SO(2n)$, provide with the symmetric metric $\breve{g}$ represented by the inner product $(\cdot,\cdot)=-B|_{\mathfrak{q}}$, $\mathfrak{q}$ horizontal space of the original fibration. According \ref{specBn}, the spectrum of the Laplacian $\Delta_{\breve{g}}$ on $(S^{2n},\breve{g})$ is $$ \sigma(\Delta_{\breve{g}})=\{\beta_q=\dfrac{q(q+2n-1)}{2(2n-1)};q\in\mathbb{N}\}.$$
 
The degeneracy instants for $\cvh$ in $]0,b]$ are the real values $t^{\text{\bf h}}_{q}$, solutions of the equation $$\dfrac{\textrm{scal}(t)}{m-1}=\beta_q=\dfrac{q(q+2n-1)}{2(2n-1)}. $$ The explicit formula for $t^g_q$, which represents the solution of the above equation in $t$, is presented in \eqref{seqso1}. Note that the constant eigenvalues $\beta_q=\dfrac{q(q+2n-1)}{2(2n-1)}$ tend to $+\infty$ as $q\rightarrow \infty$ and, since $\dfrac{\text{scal}(t)}{m-1}$ is continuous and tends to $+\infty$ as $t\rightarrow 0$, $t^{\text{\bf h}}_{q}\rightarrow 0.$ \cqd 

By a totally analogous argument, we also have the Morse index $N(t)$ for each $\cvh$, as we obtained for $\cvg$ above.
\begin{proposition}\label{morseindexso}
The Morse index of $\cvh$ is given by 
\begin{equation} 
N(t)=\left\{
\begin{array}{rcl}
\displaystyle\Sigma_{q=1}^r m_q(S^{2n}),& \mbox{if}&t^{\text{\bf g}}_{r+1}\leq t < t^{\text{\bf g}}_{r}\\ 
0, &\mbox{if}& t^{\text{\bf g}}_1\leq t\leq 1
\end{array}
\right.,
\end{equation}where $m_q(S^{2n})$ is the multiplicity of the $q$th eigenvalue of the basis $S^{2n}=SO(2n+1)/SO(2n).$ 
\end{proposition}

As consequence of the above results, we have obtained the bifurcation and rigidity instants for \linebreak$(SU(n+1)/T^{n},\cvg)$ and $(SO(2n+1)/T^n,\cvh)$, $0<t\leq 1,$ in the following theorem.

\begin{theorem} For the canonical variations $(SU(n+1)/T^{n},\cvg)$ and $(SO(2n+1)/T^n,\cvh)$, constructed above, the elements of the sequences $\{t^{\text{\bf g}}_{q}\},\{t^{\text{\bf h}}_{q}\}\subset \, ]0,b]$, of degeneracy instants for $\cvg$ and $\cvh$, given in \eqref{seqsu} and \eqref{seqso1}, are bifurcation instants for $\cvg$ and $\cvh$, respectively. Moreover, $\cvg$ and $\cvh$ are locally rigid for all $t\in]0,1]\: \setminus \: \{t^{\text{\bf g}}_{q}\}$(respectively $ t\in]0,1]\: \setminus \: \{t^{\text{\bf h}}_{q}\})$.
\end{theorem}

{\bf Proof:} It is known that, if $t_{\ast}>0$ is not a degeneracy instant for $\cvg$ and $\cvh$, then $\cvg$ and $\cvh$ are locally rigidity at $t_{\ast}.$ In the previous theorems, we proved that the degeneracy instants for $\cvg$ and $\cvh$ form the sequences $t^{\text{\bf g}}_{q},t^{\text{\bf h}}_{q}$, and thus, if $t\notin \{t^{\text{\bf g}}_{q}\}$ or $t\notin \{t^{\text{\bf h}}_{q}\}$ , $t$ must be a local rigidity instant for $\cvg$ and $\cvh$, respectively. 

The fact that each $t^{\text{\bf g}}_q$ or $t^{\text{\bf h}}_q$ are a bifurcation instants follows from Proposition \ref{propbifurcation} and Lemma \ref{autoconstante}. \cqd

The eigenvalues of the Laplacian on the basis spaces $Sp(n)/U(n)$ and $SO(2n)/U(n)$ are equal to a polynomial depending on several (integers) variables, so that we can not obtain general formula for bifurcation instants in these cases, as for the canonical variations $\cvg$ and $\cvh$ defined on $SU(n+1)/T^{n}$ and $SO(2n+1)/T^n$, respectively.

\begin{theorem} For the canonical variations $(Sp(n)/T^n,\cvk), n\geq 5,$ and $(SO(2n)/T^n,\cvm), \linebreak n\geq 4,$ introduced above, the bifurcation instants form discrete sets $\{t^g_{x_1x_2\ldots x_l}\}_{x_1,x_2,\ldots , x_l\in\mathbb{Z}_{+}}\subset \, ]0,1]$,  $1\leq l\leq n$, with infinite elements accumulating close to zero as $x_1,x_2,\ldots, x_l$ vary over $\mathbb{Z}_{+}$. Moreover, $\cvk$ and $\cvm$ are locally rigid for all $t\in]0,1]\: \setminus \: \{t^g_{x_1x_2\ldots x_l}\}_{x_1,x_2,\ldots , x_l\in\mathbb{Z}_{+}}$.
\end{theorem}

{\bf Proof:} We must find the instants $0<t<1$ such that $\dfrac{\textrm{scal}(t)}{m-1}$ is a constant eigenvalue, where $m$ represents the dimension of each spaces $Sp(n)/T^n$ and $SO(2n)/T^n$, while $\textrm{scal}(t)$ denotes the scalar curvature of $\cvk$ and $\cvm$. 

We have that $\lambda^{k,0}(t)$ is eigenvalue of $\Delta_t$ if and only if $\lambda^{k,0}(t)$ belongs to the spectrum of the Laplacian on the basis space provide with the symmetric metric $\breve{g}$ represented by the inner product $(\cdot,\cdot)=-B|_{\mathfrak{q}}$, $\mathfrak{q}$ horizontal space of the original fibration. The spectrum of the Laplacian $\Delta_{\breve{g}}$ on the basis is given by $$\sigma(\Delta_{\breve{g}})=\{\mu(\Lambda)=-B(\Lambda+2\delta,\Lambda);\Lambda\in D(G,H)\},$$ implying that any eigenvalue of the Laplacian on basis spaces $Sp(n)/U(n)$ and $SO(2n)/U(n)$ is equal to a polynomial expression depending on several (integers) variables, as has been seen in the above remark. Let $\mu(\Lambda)=\mu(x_1,\ldots,x_l)\in\sigma(\Delta_{\breve{g}})$, for some $1\leq l\leq n,$ a constant eigenvalue of the Laplacian $\Delta_t$ acting on  $(Sp(n)/T^n,\cvk), n\geq 5,$ or $(SO(2n)/T^n,\cvm), n\geq 4.$
 
The degeneracy instants for $\cvk$ and $\cvm$ in $]0,1[$ are the real values $t^g_{x_1x_2\ldots x_l}$ which are solutions of the equation $$\dfrac{\textrm{scal}(t)}{m-1}=\mu(x_1,\ldots,x_l).$$The above equations in $t$, both for $\cvk$ and $\cvm$, have infinite solutions, since the functions $]0,1[\ni t\mapsto \dfrac{\textrm{scal}(t)}{m-1} $ are continuous and assume any positive value greater than or equal to the first positive eigenvalue of the Laplacian $\Delta_{\breve{g}}$ on the respective basis spaces $Sp(n)/U(n)$ and $SO(2n)/U(n)$, namely, greater than 1. 

The fact that the degeneracy instants $t^g_{x_1x_2\ldots x_l}$ are bifurcation instants follows from Proposition \ref{propbifurcation} and Lemma \ref{autoconstante}.

We know that, if $t_{\ast}>0$ is not a degeneracy instant for $\cvk$ and $\cvm$, then $t_{\ast}>0$ is a local rigidity instant for them. Therefore, $\cvk$ and $\cvm$ are locally rigid at all $$t\in]0,1]\: \setminus \: \{t^g_{x_1x_2\ldots x_l}\}_{x_1,x_2,\ldots , x_l\in\mathbb{Z}_{+}}.$$ 

Note that in fact $t^g_{x_1x_2\ldots x_l}\in ]0,1[$, since we obtain positive numbers $0<b<1$ in Lemma \ref{autoconstante} such that $\cvk$ and $\cvm$ are locally rigid at all $t\in \, ]b,1]$. Thus, $$t^g_{x_1x_2\ldots x_l}\in\, ]0,b[ \subset ]0,1[.$$\cqd
  
It was established that the homogeneous fibration $$\pi:(G_2/T,g)\longrightarrow (G_2/SO(4),\breve{g})$$ is a Riemannian submersion with totally geodesic fibers isometric to $(H/K,\hat{g})$, $H/K=SO(4)/T\cong SO(4)/SO(2)\times SO(2)\cong S^2\times S^2$, $g$ the normal homogeneous metric determined by the inner product $(-B)|_{\m}$, $\hat{g}$ the metric given by $(-B)|_{\mathfrak{p}}$ and $\breve{g}$ defined by the inner product $(-B)|_{\mathfrak{q}}$. Furthermore, the canonical variation $g_t$ of this submersion is represented by the inner product {\small $$(\cvn)_{eT}=(-B)|_{\m_{\alpha_1}}+(-B)|_{\m_{\alpha_2}}+t^2(-B)|_{\m_{\alpha_1+\alpha_2}}+(-B)|_{\m_{\alpha_1+2\alpha_2}}+t^2(-B)|_{\m_{\alpha_1+3\alpha_2}}+(-B)|_{\m_{2\alpha_1+3\alpha_2}}.$$} 

The first positive eigenvalue of the Laplacian $\Delta_{\hat{g}}$ on the fiber is, according to \cite{Urakawa}, equal to $\phi_1=1$ and the first positive eigenvalue of the Laplacian $\Delta_g$ on the total space is equal to $\mu_1=\dfrac{1}{2}$. Hence, $$\lambda^{1,1}(t)=\mu_1+(\frac{1}{t^2}-1)\phi_1=\dfrac{1}{2}+(\frac{1}{t^2}-1).$$

In Proposition \ref{propscalG}, we obtained the expression \eqref{scalG} of $\text{scal}(t)$, scalar curvature of the canonical variation $\cvn$, namely  $$\text{scal}(t)=\dfrac{2+12t^2-2t^4}{3t^2}, \ \ t>0.$$ Since $m=\dim G_2/T=12$, it follows that $$\dfrac{\text{scal}(t)}{m-1}=\dfrac{2+12t^2-2t^4}{33t^2}.$$ Furthermore, we have the spectrum of the Laplacian on the basis space $(G_2/SO(4),\breve{g})$  $$\sigma(\Delta_{\breve{g}})=\left\{\dfrac{1}{6}(9r+6r^2+5s+6rs+2s^2);\mathbb{Z}\ni r,s\geq 0\right\}.$$

Now, we can determine the bifurcation and local rigidity instants for $\cvn$ in the interval $]0,1]$.

\begin{theorem} The elements of the set $\{t^{\text{\bf n}}_{rs}\}\subset ]0,1]$ are given by {\small $$t^{\text{\bf n}}_{rs}=\dfrac{\sqrt{\sqrt{\left(-66 r^2-33 r s-99 r-22 s^2-55 s+24\right)^2+64}-66 r^2-33 r s-99 r-22 s^2-55 s+24}}{2 \sqrt{2}},$$}$\mathbb{Z}\ni r,s\geq 0$, are bifurcation instants for $(G_2/T,\cvn)$.  Moreover, $\cvn$ is locally rigid for each $0<t\leq 1$ such that $t\notin \{t^{\text{\bf n}}_{rs}\}$.
\end{theorem}

{\bf Proof:} By Proposition \ref{propbifurcation}, if $0<t_{\ast}<1$ such that $\dfrac{\text{scal}(t_{\ast})}{m-1}\in\sigma(\Delta_{\breve{g}})$ and $\dfrac{\text{scal}(t_{\ast})}{m-1}<\lambda^{1,1}(t_{\ast})$, then $t_{\ast}$ is a bifurcation instant.

We have $$\dfrac{\text{scal}(t)}{m-1}=\dfrac{\text{scal}(t)}{11}=\dfrac{2+12t^2-2t^4}{33t^2}<\lambda^{1,1}(t)=\mu_1+(\frac{1}{t^2}-1)\phi_1=\dfrac{1}{2}+(\frac{1}{t^2}-1),$$for all $0<t\leq 1$. The instants that satisfy $\dfrac{\text{scal}(t_{\ast})}{11}\in\sigma(\Delta_{\breve{g}})\subset\sigma(\Delta_t)$ are the solutions of $$\dfrac{\text{scal}(t)}{11}=\dfrac{2+12t^2-2t^4}{33t^2}=\dfrac{1}{6}(9r+6r^2+5s+6rs+2s^2), 0<t\leq 1$$ which are exactly the $t^{\text{\bf n}}_{rs}$ given above. The elements of the set $\{t^{\text{\bf n}}_{rs}\}$ are such that $$\dfrac{\text{scal}(t^{\text{\bf n}}_{rs})}{11}<\lambda^{1,1}(t^{\text{\bf n}}_{rs})$$ and $t^{\text{\bf n}}_{rs}$ is a bifurcation instant for all $0\leq r,s\in\mathbb{Z}.$ It follows that $\cvn$ is locally rigid at all $0<t<1$ such that $t\notin \{t^{\text{\bf n}}_{rs}\}$, i.e, at every instant in the complement of the bifurcation instants. \cqd

\begin{remark} By continuity of $\dfrac{\text{scal}(t)}{11}$ and by the fact that $\lim_{t\rightarrow 0}\dfrac{\text{scal}(t)}{11}=+\infty$, since the eigenvalues of $\Delta_{\breve{g}}$ goes to $+\infty$ when $r,s\rightarrow +\infty$, we have that the sequence $t_{rs}^g$ obtained above is such that $t^{\emph{\bf n}}_{rs}\rightarrow 0$ when $r,s\rightarrow +\infty$ and, then, we determine a sequence of bifurcation instants accumulating close to zero. 
%\begin{figure}[!h]
%\centering
%\includegraphics[height=5cm,width=7cm]{excepcional.jpg}
%\caption{Graphic of $\dfrac{\text{scal(t)}}{11}$}
%\end{figure}
\end{remark}
%\begin{figure}[h]
%\centering
%\includegraphics[height=7cm,width=10cm]{G2.jpeg}
%\caption{For $G_2/T$, graphs of the functions $\dfrac{\textrm{scal}(t)}{m-1}$ in red and $\lambda^{k,j}(t)$ in black, the constants corresponding to $(k,0), 1\leq k\leq 6$, and non-constants to $(k,j)$ where $1\leq j\leq k\leq 6$. The dashed vertical lines mark five degeneracy instants (which are all bifurcation instants) starting at $t^{\text{\bf n}}_{1,0}$.}
%\end{figure}

\section{Multiplicity of Solutions to the Yamabe Problem}\label{Section3.3}

We now explain how to obtain multiplicity results for the canonical variations $\cvg,\cvh,\cvk,\cvm$ and $\cvn$ applying the next proposition due Bettiol and Piccione in \cite{pacificjournal}. We have been interested in determine which conformal classes carry multiple unit volume metrics with constant scalar curvature. 

\begin{proposition}[\cite{pacificjournal}] Let $g_t$, with $t\in \, ]0,b[$, be a family  of metrics on $M$ with $N(t)=N(g_t)>0$, $N(t)$ the Morse index of $g_t,$ and suppose there exists a sequence $\{t_q\}$ in $]0,b[$, that converges to 0, of bifurcation values for $g_t$. Then, there is an infinite subset $\mathcal{G}\subset ]0,b[$ accumulating at $0$, such that for each $t\in\mathcal{G}$, there are at least $3$ solutions to the Yamabe problem in the conformal class $[g_t].$
\end{proposition}

Applying the last Proposition and our bifurcation results, we can determine a lower bound for the number of unit volume metrics with constant scalar curvature in each conformal class $[\cvg], [\cvh], [\cvk], [\cvm], [\cvn]$, respectively, for instants in a given subset $\mathcal{G}\subset ]0,1[$. 

\begin{theorem}
Let $\cvg,\cvh,\cvk,\cvm$ and $\cvn$ be the families of homogeneous metrics described above. Then, there exists, for each of such family, a subset $\mathcal{G}\subset ]0,1[$, accumulating at $0$, such that for each $t\in\mathcal{G}$, there are at least $3$ solutions to the Yamabe problem in each conformal class $[\cvg], [\cvh], [\cvk], [\cvm], [\cvn].$
\end{theorem}

{\bf Proof:} It is only necessary to verify that there exists $0<b<1$ such that $N(t)>0$ in $]0,b[$. For each canonical variation $\cvg,\cvh,\cvk,\cvm$ and $\cvn$ take a positive number $0<b<1$ such that $\dfrac{\text{scal(t)}}{m-1}<\lambda_1(t)$ for all $t\in\, ]b,1[$, with $m$ denoting the dimension of the respective total space. This implies that $N(t)=0$, for all $t\in\, ]b,1[$, since there are no eigenvalues of $\Delta_t$ less than  
$\dfrac{\text{scal(t)}}{m-1}$. If $t=b$, one has $\dfrac{\text{scal(t)}}{m-1}=\beta_1,$ where $\beta_1$ is the first positive eigenvalue of the Laplacian on the basis. We proved that $t=b$ is a bifurcation instant and then the Morse index changes from $0$ to a positive integer. Hence, for $t\in\, ]0,b[\subset]0,1[,$ we have $N(t)\geq N(b-\epsilon)>0$, since by definition, the Morse index is equal to the number (counting multiplicity) of positive eigenvalues that are lees than $\dfrac{\text{scal}(t)}{m-1}$, which is strictly decreasing for $0<t<1$ and $\dfrac{\text{scal}(t)}{m-1}\rightarrow \infty $ as $t\rightarrow 0.$, for $\cvg,\cvh,\cvk,\cvm$ and $\cvn$. \cqd
\bibliography
{unsrt}  


\end{document}